\documentstyle{amsart}

\makeatletter
\def\verbatim{\interlinepenalty\@M \@verbatim
  \leftskip\@totalleftmargin\advance\leftskip2pc
  \frenchspacing\@vobeyspaces \@xverbatim}
\makeatother
\newtheorem{theorem}{Theorem}[section]
\newtheorem{corollary}[theorem]{Corollary}
\newtheorem{lemma}[theorem]{Lemma}
\newtheorem{prop}[theorem]{Proposition}

\theoremstyle{definition}
\newtheorem{definition}{Definition}[section]

\numberwithin{equation}{section}

\newcommand{\clp}{{\rm cl}_{{\Bbb R}^p} \hspace{.02in}}

\newcommand{\ra}{\rightarrow}

\newcommand{\va}{\left|}
\newcommand{\vb}{\right|}
\newcommand{\vc}{\left\|}
\newcommand{\vd}{\right\|}

\newcommand{\hs}{\hspace{.02in}}

\newcommand{\noe}{A^n (\Omega, E)}
\newcommand{\nmoe}{A^n (\Omega', F)}

\newcommand{\coe}{C^n (\Omega, E)}
\newcommand{\no}{C^n (\Omega, {\Bbb  R})}
\newcommand{\rno}{A^n (\Omega, {\Bbb  R})}

\newcommand{\noor}{C^n_c (\Omega, {\Bbb  R})}

\newcommand{\moe}{A^m (\Omega', F)}
\newcommand{\bnoe}{C^n_* (\Omega, E)}

\newcommand{\bmoe}{C^m_* (\Omega', F)}

\title[Automatic continuity]{Automatic continuity and weighted composition operators between spaces of vector-valued differentiable functions}

\author{Jes\'us Araujo}
\address{Departamento de Matem\'aticas,
Estad\'{\i}stica y Computaci\'on\\ Universidad de Cantabria\\
Facultad de Ciencias\\ Avda.
de los Castros, s. n.\\ E-39071 Santander, Spain}

\email{araujoj@@unican.es}
\thanks{Research partially supported by
the Spanish Direcci\'on General de Investigaci\'on Cient\'{\i}fica
y T\'ecnica (DGICYT, PB98-1102).}

\date{}                                          
\begin{document}

\maketitle

\begin{abstract}
Let $E$ and $F$ be Banach spaces.
It is proved that if $\Omega$
and $\Omega'$ are open subsets of ${\Bbb R}^p$ and ${\Bbb R}^q$, respectively, and $T$ is a linear biseparating map between two spaces of differentiable functions $\noe$ and $\moe$,
then $p = q$, $n=m$, and there exist a diffeomorphism $h$ of
class $C^m$ from $\Omega'$ onto $\Omega$, and a map $J: \Omega' \ra L(E,F)$ of class ${\rm s}- C^n$ such that for every $y \in \Omega'$ and every $f \in
\noe$, $(Tf) (y) = (Jy) (f(h(y)))$. In particular $E$ and $F$ are isomorphic as Banach spaces and, as a consequence, all linear biseparating maps are continuous for usual topologies in the spaces of differentiable functions. 
\end{abstract}

\section{Introduction}

It is well known that an algebraic link between spaces of continuous functions may lead to a topological link between the spaces on which the functions are defined. For instance, it turns out that if there exists a ring isomorphism $T :C(X) \ra C(Y)$, then the realcompactifications of $X$ and $Y$ are homeomorphic (\cite[pp. 115-118]{GJ}). Also if $h$ is the resultant homeomorphism from the realcompactification of $Y$ onto that of $X$, then $Tf = f \circ h$ for every $f \in C(X)$, so we have a complete description of it. As a result, when $X$ and $Y$ are realcompact, we deduce that if both spaces of continuous functions $C(X)$ and $C(Y)$ are endowed with the compact-open topology, then every ring isomorphism between them is continuous. In this  result, the key point is that every ring isomorphism sends maximal ideals into maximal ideals. This implies that a good description of maximal ideals lead to the definition of a map from $Y$ onto $X$. 

Of course, the pattern above have been succesfully applied to many other algebras of functions. However the situation becomes more complicated if we consider spaces of functions which take values in {\em arbitrary} Banach spaces. In this context and unlike algebra or ring homomorphisms, we can still use mappings satisfying the property $\vc Tf \vd \vc Tg \vd \equiv 0$ if and only if $\vc f \vd \vc g \vd \equiv 0$. These maps are called {\em biseparating}, and coincide with disjointness preserving mappings whose inverses preserve disjointness too (\cite{AK}). In general these maps turn out to be efficacious substitutes for homomorphisms. Indeed, in \cite{A1}, we prove that the existence of a biseparating mapping between a large class of spaces of vector-valued continuous functions $A(X,E)$ and $A(Y,F)$ ($E$, $F$ are Banach spaces) yields homeomorphisms between some compactifications (and even the realcompactifications) of $X$ and $Y$. The automatic continuity of a {\em linear} biseparating mapping is also accomplished in some cases (see \cite{A2,A3}). Related results have also been given recently, for some other families of scalar-valued functions, for instance, in \cite{ABN, FH, J} and \cite{JW}. In this paper, we go a step beyond and work in a context which does not seem to have made its way into the literature yet, namely, linear operators between spaces of differentiable functions taking values in arbitrary Banach spaces.

As for spaces (indeed algebras) of {\em scalar-valued} differentiable functions, Myers (\cite{My}) showed that the structure  of a compact differentiable manifold of class $C^n$ is determined by the algebra of all real-valued functions on $M$ of class $C^n$. In this line, Pursell (\cite{P}) checked that the ring structure of infinitely differentiable functions defined on an open convex set of ${\Bbb R}^n$ determines such set up to a diffeomorphism. Homomorphisms between algebras of differentiable functions defined on real Banach spaces have been studied by Aron, G\'omez and Llavona (\cite{AGL}); in this paper a description of homomorphisms is given and the automatic continuity is obtained as a corollary in quite a general setting, which includes in particular the case when the Banach spaces are finite-dimensional. Further,  Guti\'errez and Llavona (\cite{GL}) pursued the study  of continuous homomorphisms between algebras of differentiable functions defined on open subsets of real Banach spaces, obtaining a description in many cases. On the other hand, automatic continuity results for algebras of differentiable functions have been also given for instance in \cite{BCL}, \cite{Lr}, \cite{KN} and \cite{NRV}. For classical results and techniques in the study of automatic continuity, see also \cite{D} and \cite{S}. Finally, Molnar states in the final remarks of \cite{Ml} that it would seem to be of some importance to study some of the results above for separating mappings instead of algebra or ring isomorphism. In his opinion, that is certainly more difficult to carry out.

Our main goal here is to investigate the behaviour of biseparating maps when defined between spaces of vector-valued differentiable functions. Let $E$ and $F$ be Banach spaces. We shall prove that if $\Omega$ and $\Omega'$ are open subsets of ${\Bbb R}^p$ and ${\Bbb R}^q$, respectively, and $T$ is a linear biseparating map between two spaces of differentiable functions $\noe$ and $\moe$, then $p = q$, $n=m$, and there exist a diffeomorphism $h$ of class $C^m$ from $\Omega'$ onto $\Omega$, and a map $J: \Omega' \ra L(E,F)$ of class ${\rm s}- C^n$ such that for every $y \in \Omega'$ and every $f \in \noe$, $(Tf) (y) = (Jy) (f(h(y)))$. In particular $E$ and $F$ are isomorphic as Banach spaces and, as a consequence, all linear biseparating maps are continuous for usual topologies in the spaces of differentiable functions.

\medskip

        In Section 2 (as well as in the rest of Section 1), we assume that $E$ is a real Banach space.

 If $\lambda = (\lambda_1 , \lambda_2 ,
\ldots , \lambda_p)$ is a $p$-tuple of non-negative integers, we
set $\va \lambda \vb := \lambda_1 + \lambda_2 + \ldots +
\lambda_p$. If $\Omega$ is an open subset of ${\Bbb R}^p$, then
$\coe$ consists of the $E$-valued functions $f $ in $\Omega$
whose partial derivatives \[{\partial^{\lambda
} f} := \frac{\partial^{\lambda_1 +
\lambda_2 + \ldots + \lambda_p}
f}{\partial x_1^{\lambda_1} \partial x_2^{\lambda_2} \ldots
\partial x_p^{\lambda_p}}\]exist and are continuous for each $\lambda =   (\lambda_1
, \lambda_2 , \ldots , \lambda_p) \in \Lambda$,
where $\Lambda
:= \{ \lambda \in ({\Bbb N} \cup \{0\})^p : \va \lambda \vb 
\le n \}$. It is well known that, if for $L^k ({\Bbb R}^n , E)$
we denote the space of continuous $k$-${\Bbb R}$-linear maps of ${\Bbb R}^n$
into $E$, then $\coe$ coincides with the space
of maps $f: \Omega \ra E$ such that the differential $D^k f:
\Omega \ra L^k ({\Bbb R}^n , E)$ exists and is continuous for
each $k=0, \ldots, n$.

        From now on we will assume that $n \ge 1$, and that $\Omega \subset {\Bbb R}^p$ is a (nonempty) open set. Also, for $C \subset {\Bbb R}^p$, $\clp C$ denotes its closure in ${\Bbb R}^p$. $\noor$ will denote the subring of
$\no$, respectively, of functions with compact
support. On the other hand,  if $x = (x_1 , x_2 , \ldots, x_p)$ belongs to
${\Bbb R}^p$, we set $\va x \vb := \max_{j} \va x_j \vb$.

\section{Some previous results}

\begin{lemma}\label{bueno}
Suppose that $a_0 , a_1 , \ldots, a_k \in \coe$, and that $f:
\Omega \times {\Bbb R} \ra E$ is a polynomial in $t$
defined as\[f(x, t) = 
\sum_{i=0}^k a_i (x) t^i\]for every $(x,t) \in \Omega \times
{\Bbb R}$. Then $f \in C^n (\Omega \times {\Bbb R}, E)$.
\end{lemma}

\begin{pf}
It is immediate from the fact that all partial derivatives up to
order $n$ exist and are continuous.
\end{pf}

\begin{lemma}\label{dabute}
Suppose that $f, g \in C^n (\Omega, E)$ and $k : \Omega \ra {\Bbb R}$
satisfy \[f(x) = k(x) g(x)\] for every $x \in \Omega$. If
$g(x) \neq 0$ for every $x \in \Omega$, then $k \in C^n (\Omega, {\Bbb R})$.
\end{lemma}

\begin{pf}
Fix $x_0 \in \Omega$. We are going to prove that $k$ is of
class $C^n$ in a neighborhood of $x_0$. First we take ${\bf f}' \in E'$
such that ${\bf f}' (g (x_0)) \neq 0$. Since
${\bf f}'$ and $g$ are continuous, then there exists an
open neighborhood $U$ of $x_0$ such that
\[{\bf f}' (g (x)) \neq 0\] for every $x \in U$.

        Now
\[{\bf f}' (f (x)) = k(x) {\bf f}' (g (x))\] for every $x \in
U$. This implies that, for every $x \in U$, 
\[k (x) = \frac{{\bf f}' (f(x))}{{\bf f}'
(g(x))},\] which is the quotient of two real-valued
functions of class $C^n$.

        This proves that $k$ is of class $C^n$.
\end{pf}

It is well known that, when it exists at a point $a \in \Omega$,
the differential $D^k f (a) \in L^{k} ({\Bbb R}^p , E)$ is a
symmetric form of degree
$k$.

Now, given a map $f \in \coe$, we define its Taylor polynomial
function of degree $n$ at $a \in \Omega$ as 
\begin{eqnarray*}
T_a (x) &:=& f(a) +
Df(a) (x-a) + \frac{1}{2} D^2 f(a) (x-a, x-a) + \ldots \\&+&
\frac{1}{n!} D^n f(a) (\underbrace{x-a, x-a, \ldots, x-a}_{n}).
\end{eqnarray*}

The proof of the following result is straightforward.

\begin{lemma}\label{chino}
Suppose that $f \in \coe$, $a \in \Omega$, and $1 \le k \le n$.
Then the $k$-th derivative of the Taylor polynomial function
$T_a$ of
degree $n$ of $f$ is equal to the Taylor polynomial
function of degree $n-k$ of $D^k f$ at $a$.
\end{lemma}

The following theorem, known as {\em Whitney's extension theorem}, can
be found, for instance, in \cite[Theorem 3.1.14]{Fr}.

\begin{theorem}\label{47}
Suppose $n \in {\Bbb N}$, $A$ is a
closed subset of ${\Bbb R}^p$, and to each $a \in A$ corresponds
a polynomial function \[ P_a : {\Bbb R}^p \ra E\] with degree
$P_a \le n$. Whenever $C \subset A$ and $\delta >0$ let $\rho
(C, \delta)$ be the supremum of the set of all numbers\[ \vc D^i
P_a (b) - D^i P_b (b) \vd \cdot \va a - b \vb^{i-n} \cdot
(n-i)!\]corresponding to $i = 0, \ldots , n$ and $a, b \in C$
with $0 < \va a-b \vb \le \delta$.

        If $\rho (C, \delta) \ra 0$ as $\delta \ra 0+$ for each
compact subset $C$ of $A$, then there exists a map $g: {\Bbb R}^p
\ra E$ of class $C^n$ such that \[D^i g(a) = D^i P_a (a)\] for
$i=0, \ldots , n$ and $a \in A$.
\end{theorem}

\begin{prop}\label{45}
Let $p \ge 2$. For $ s \in {\Bbb R}$, $s>0$, consider the following
compact subsets of ${\Bbb R}^p$:
\[A := \{ (x_1 , x_2 ,
\ldots , x_p ) \in {\Bbb R}^p :
\sum_{i=2}^p x_i^2 \le x_1^2 / 9 , \va x_1 \vb \le s \},\]
\[A^+ := \{ (x_1 , x_2 ,
\ldots , x_p ) \in A : x_1 \ge 0\},
\]
and \[ A^- :=
\{ (x_1 , x_2 ,
\ldots , x_p ) \in A : x_1 \le 0\} .\]

Suppose that $\Omega$ is an open subset of ${\Bbb R}^p$
containing $A$ and
that $f$ belongs to $\coe$.
If  \[\partial^{\lambda} f  (0, 0, \ldots , 0) =0\] for every
$\lambda \in \Lambda$, then there
exists a function $f_+ \in \coe$ with compact support such that, for $\lambda \in \Lambda$, \[\partial^{\lambda} f_+  (x) = \partial^{\lambda} f (x) \]for every $x \in A^+$, and
\[\partial^{\lambda} f_+   (x) = 0\] for every
$x \in A^-$.
\end{prop}

\begin{pf}
Suppose that for $x \in A$,
$T_x$ stands for the polynomial function of degree $n$ given in the
Taylor formula for $f$ at $x$. Now for $a \in A^-$, we consider
as $P_a$ the polynomial identically zero, and
for $a \in A^+$, we consider as
$P_a$ the polynomial $T_a$. As it is seen for instance in \cite[Theorem 2.71]{Ch} or
\cite[p. 350]{Ln}, if $r>0$ and $\va b-a \vb < r$, we have that
\[ \vc T_a (b) - f (b) \vd \le \frac{\va b-a \vb^{n}}{ n!}
\sup_{\va x-a \vb \le r} \vc D^n f(x) - D^n f(a) \vd.\]Now it is
easy to see that, by Lemma~\ref{chino},  for $i
\in \{1, 2, \ldots , n\}$,
\[ \vc D^i
T_a (b) - D^i f (b) \vd \le \frac{\va b-a \vb^{n-i}}{ (n-i)!}
\sup_{\va x-a \vb \le r} \vc D^n f(x) - D^n f(a) \vd .\]

        This proves that if $\va b-a \vb <r$, 
\[ \vc D^i
T_a (b) - D^i f (b) \vd \cdot \va b - a \vb^{i-n} \cdot
(n-i)! \le  \sup_{\va x-a \vb \le r} \vc D^n f(x) - D^n f(a)
\vd \] for $i
\in \{0, 1, \ldots , n\}$. Then we  define $\tau (r)$ as the
supremum of the set
of all numbers \[   \sup_{\va x-a \vb \le r} \vc D^n f(x) - D^n f(a)
\vd \]for $x, a \in A$, which is a real number, because $A$ is compact.
Clearly, since $D^n f$ is continuous, if $r$ tends to zero,
$\tau(r)$ tends to zero. Now
suppose that $a$ and $b$ belong to $A$ and $0 < \va b -a \vb \le
r$. Then we have the following possibilities:

\begin{itemize}
\item $a , b \in A^+$. Then
we have that, by Lemma~\ref{chino}, for $i \in \{0, 1, \ldots , n\},$ \[ D^i P_b (b) =
D^i f(b) \] and
consequently \[ \vc D^i P_a (b) - D^i P_b (b) \vd = \vc D^i T_a
(b) - D^i f (b) \vd .\]
\item $a, b \notin A^+$. Then \[ \vc D^i P_a (b) - D^i P_b (b) \vd =0 .\]
\item $a \notin A^+$, $b \in A^+$. Note that since we are assuming by hypothesis that $D^i f(0, 0,
\ldots, 0) =0$,
then $D^i P_0 (b) =0$
for all $i \in \{0, 1, \ldots , n\}$. Consequently
\begin{eqnarray*}
\vc D^i P_a (b) - D^i
P_b (b) \vd &=& \vc D^i P_b (b) \vd \\ &=& \vc D^i P_0 (b) - D^i P_b
(b) \vd \\ &=& \vc D^i T_0 (b) - D^i f(b) \vd.
\end{eqnarray*}
\item $a \in A^+$, $b \notin A^+$. Then we have that
\begin{eqnarray*}
\vc D^i P_a (b) - D^i
P_b (b) \vd &=& \vc D^i T_a (b) \vd \\ &\le& \vc D^i T_a (b) - D^i f
(b) \vd + \vc D^i f (b) \vd \\  
&=& \vc D^i T_a (b) - D^i f
(b) \vd + \vc D^i f (b) - D^i T_0 (b) \vd.
\end{eqnarray*}
\end{itemize}

        Note that in the third and forth cases above, $\va a
\vb, \va b \vb \le \va b -a \vb \le r$. This implies that in
these two cases
\begin{eqnarray*}
\vc D^i
P_a (b) - D^i P_b (b) \vd \cdot \va b - a \vb^{i-n} \cdot
(n-i)! &\le&  \sup_{\va x-a \vb \le r} \vc D^n f(x) - D^n f(a)
\vd \\ &+& \sup_{\va x \vb \le r} \vc D^n f(x) - D^n f(0)
\vd.
\end{eqnarray*}

        On the other hand it is easy to see that in the other
two cases \[ \vc D^i
P_a (b) - D^i P_b (b) \vd \cdot \va b - a \vb^{i-n} \cdot
(n-i)! \le  \sup_{\va x-a \vb \le r} \vc D^n f(x) - D^n f(a)
\vd.\] 

This facts imply that, if $\rho$ is defined as in Theorem~\ref{47}, $\rho (A, r) \le 2 \tau (r)$. Also, it
is clear that if $C $ is a compact subset of $A$, then $\rho (C,
r) \le \rho (A,r)$. Consequently by Theorem~\ref{47}, we have
that there exists  $f_0 \in C^n (\Omega, E)$ such that, given any $\lambda \in \Lambda$,
\[\partial^{\lambda} f_0  (x) = \partial^{\lambda} f   (x) \]for every $x \in
A^+$, and \[\partial^{\lambda} f_0  (x) = 0\] for every
$x \in A^-$. Also it is clear that if we take $g_0 \in \noor$ such that $g_0 \equiv 1$ on an open neighborhood of $A$, then $f_+ := g_0 f_0 \in \coe$ satisfies the requirements of the theorem.
\end{pf}

\section{Biseparating maps: a first approach}

In previous sections, we assumed that $E$  was a {\em real} Banach space. Imagine now that it is a {\em complex} Banach space instead. It is clear that in any case, it can also be viewed as a {\em real} space, and in this sense we consider defined the space $\noe$. It is immediate that all results given before hold for them. 

But being also complex ensures that  $\noe$ is both real and complex as a linear space, and consequently we can consider both real and complex linear maps from $\noe$ into some  other vector spaces. This is the reason why, unlike so far, 
in this section we assume that
$E$ and $F$ are ${\Bbb K}$-Banach spaces, where ${\Bbb K} =
{\Bbb R}$ or ${\Bbb C}$.

\medskip

In this section (and in the following ones) we also assume that our spaces are defined according to one of the following two situations.

\begin{itemize}
\item {\bf Situation 1.} $\Omega$ and $\Omega'$ are (not necessarily bounded) open subsets of ${\Bbb R}^p$ and ${\Bbb R}^q$, respectively ($p, q \in {\Bbb N}$). In this case,  $\noe = \coe$ and $\moe = C^m (\Omega', F)$ ($n, m \ge 1$).
\item {\bf Situation 2.}  $\Omega$ and $\Omega'$ are {\em bounded} open subsets of
${\Bbb R}^p$ and ${\Bbb R}^q$, respectively ($p, q \in {\Bbb N}$). In this case, $\noe = C^n
(\bar{\Omega}, E)$ and $\moe = C^m (\bar{\Omega}', F)$ ($n, m \ge 1$). Here $C^n
(\bar{\Omega}, E)$ (respectively $C^m (\bar{\Omega}', F)$) denotes the
subspace of $\coe$ (respectively $C^m (\Omega', F)$) of those functions
whose partial derivatives up to order $n$ (respectively $m$)
admit continuous extension to the boundary of $\Omega$
(respectively $\Omega'$).
\end{itemize}

\medskip

For a function $f \in \noe$ (respectively, in $\moe$), we denote by $c(f)$ the cozero set of $f$, that is, the set $\{x \in \Omega: f(x) \neq 0\}$.

As for the spaces of linear functions, we will denote by $L' (E,F)$ and by $B' (E,F)$ the sets of (not necessarily continuous) linear maps and  bijective linear maps  from $E$ into $F$, respectively. $L(E,F)$ and $B(E,F)$ will denote the spaces of {\em continuous} linear maps and  bijective {\em continuous} linear maps from $E$ into $F$.

\begin{definition}
{\em A map $T: \noe \ra \moe$ is said to be {\em separating} if it is additive and
$c(Tf) \cap c(Tg) = \emptyset$ whenever $f, g \in \noe$ satisfy
$c(f) \cap c(g) = \emptyset$. Besides $T$ is said to be
{\em biseparating} if it is bijective and both $T$ and $T^{-1}$ are separating.} 
\end{definition}

Equivalently, we see that an additive map $T: \noe \ra \moe$ is separating if $\vc (Tf) (y) \vd \vc (Tg) (y) \vd =0$ for all $y \in \Omega'$ whenever $f, g \in \noe$ satisfy $\vc f(x) \vd \vc g(x) \vd=0$ for all $x \in \Omega$. 

Let $\Omega_1 := \Omega$, $\Omega'_1 := \Omega'$ when we are in Situation 1, and $\Omega_1 := \bar{\Omega}$, $\Omega'_1 := \bar{\Omega}'$ if we are in Situation 2.  A point $x \in
\Omega_1$ is said to be a
support point of $y \in \Omega'_1$ if, for every neighborhood $U$
of $x$ in $\Omega_1$, there exists $f \in \noe$ satisfying $c(f) \subset U$
such that $ (Tf) (y) \neq 0$.

In a much more general context, (not necessarily linear) biseparating maps are studied in \cite{A1}. Some general results are given concerning conditions which allow to link certain topological spaces when dealing with a biseparating map. Suppose in particular that $A$ is a subring of the space $C(\Omega_1, {\Bbb R})$ of all real-valued functions on $\Omega_1$.  Clearly each map $f: \Omega_1 \ra {\Bbb R}$ admits a continuous extension $f^{\beta \Omega_1} : \beta \Omega_1 \ra {\Bbb R} \cup \{\infty\}$ to the Stone-\v{C}ech compactification of $\Omega_1$ (which obviously coincides with $\Omega_1$ when we are in Situation 2 above). Assume now that  the space of all extensions of elements of $A$ separate the points of $\beta \Omega_1$; following \cite{A1}, we say that  $A$  is {\em strongly regular} when given $x_0 \in \beta \Omega_1$ and a nonempty closed subset $K$ of $\beta \Omega_1$ which does not contain $x_0$, there exists $f \in A$ such that $f^{\beta \Omega_1} \equiv 1$ on a neighborhood of $x_0$ and $f^{\beta \Omega_1} \equiv 0$ on $K$. Now, let us denote by $C(\Omega_1, E)$ the space of all continuous $E$-valued functions on $E$. In particular, in \cite[Corollary 2.9]{A1}, it is proven that if $A(\Omega_1, E) \subset C(\Omega_1 , E)$ and $B (\Omega'_1, F) \subset C(\Omega'_1, F)$ are an $A$-module and a $B$-module, respectively, where $A \subset C(\Omega_1, {\Bbb R}) $ and $B \subset C(\Omega'_1 ,{\Bbb R})$ are stronly regular rings, then there exists a homeomorphism $h$ from $\Omega'_1$ onto $\Omega_1$ whenever there is a biseparating map from $A(\Omega_1, E)$ onto $B(\Omega'_1, F)$.

Applied to our situation, it is well known that if $A= C^n (\Omega_1, {\Bbb R})$ and $B= C^m (\Omega'_1, {\Bbb R})$, then $A$ and $B$ are strongly regular rings (see for instance \cite[Corollary 1.2]{W}). On the other hand, it is straightforward to see that $\noe$ and $\moe$ are an $A$-module and a $B$-module, respectively. As a consequence,  the existence of a biseparating map from $\noe$ onto $\moe$ ensures the existence of a  homeomorphism $h$  from $\Omega'_1$ onto $\Omega_1$. As it is also shown in \cite{A1}, this map $h$ sends each point in $\Omega'_1$ into its support point in $\Omega_1$, and is called {\em support map} for $T$. It turns out that the support map for $T^{-1}$ is $h^{-1}$ (\cite[Theorem 2.8]{A1}). Rephrasing Lemma 2.4 in \cite{A1}, we have the following property.

\begin{lemma}\label{txa}
If $y \in \Omega'_1$ and $f \in \noe $ vanishes on a neighborhood of $h(y)$, then $(Tf) (y) =0$.
\end{lemma}

\begin{lemma}\label{46}
Suppose that $\Omega $
contains the origin, and that $T: \noe \ra \moe$ is a  biseparating map.
Assume also that $f \in \noe$
satisfies that for all
$\lambda \in \Lambda$, \[{\partial^{\lambda
} f}  (0, 0, \ldots, 0) = 0. \]
If $(0, 0 , \ldots, 0) \in \Omega$ is the support point of
$y \in \Omega'$, 
then $(Tf) (y) =0$.
\end{lemma}

\begin{pf}
First suppose that $p>1$ and that the closed ball of center
$0$ and radius $s$ is
contained in $\Omega$. If we take $A^+$, $A^-$ and $f_+$ as in Proposition~\ref{45},
then $f_+$ and $f -f_+$ belong to $\noe$ and satisfy
$f_+ (x) = 0 $ and $(f - f_+) (x) =0$ for every $x \in A^-$ and
for
every $ x \in A^+$ respectively. We have that $(Tf) (y) = (T
f_+)(y) + (T(f - f_+))(y)$. Also,
since for any neighborhood $U$ of the origin there exists
an open subset $V$ of $U$ such that $f_+ (x)=0$ for
every $x \in V$, then we have that, taking into account that the support map for $T^{-1}$ is $h^{-1}$, by Lemma~\ref{txa}, $(T f_+ )
(h^{-1} (x)) =0$. Since $h:\Omega' \ra \Omega$ is a  homeomorphism, we deduce
that $(T f_+ ) (y) =0$ and, in the same way, $(T (f - f_+)) (y)
=0$. We conclude that $(Tf) (y) =0$.

        Consider now the case when $p=1$. Note that if $f \in \noe$
satisfies $f(0) =0$ and $0= f' (0) =  \cdots = f^{(n)} (0)$, then
it is clear that $f \xi_{(- \infty, 0)} $ and $f \xi_{(0, +
\infty)}$ belong to $\noe$ (where 
$\xi_A$ stands for
the characteristic function of $A$) and, as above, $(T(f \xi_{(0, +
\infty)})) (y) =0 = (T( f \xi_{(-\infty, 0)} ))(y)$. We conclude
that $(Tf) (y) =0$.
\end{pf}

Next lemma is a first attempt to describe all biseparating maps, and does not take into account some important details which will be discussed in Section 5.  In this way we characterize all biseparating linear maps from $\noe$ onto $\moe$ as weighted composition bijective maps. Notice that we assume no continuity properties on $T$. In fact, we will suppose that our spaces $\noe$ and $\moe$ are not endowed with any topologies. 

\begin{lemma}\label{401}
Suppose that $T: \noe \ra \moe$
is a ${\Bbb K}$-linear  biseparating map. Then $p=q$, $n=m$, and there
exist a diffeomorphism $h$ of class $C^n$ from $\Omega'$ onto
$\Omega$ and a  map
$J: \Omega' \ra B'(E,F)$ such that for every
$y \in \Omega'$ and 
every $f \in \noe$, \[ (Tf) (y) = (Jy) (f(h(y))) .\]
\end{lemma}

\begin{pf}
First, the existence of the homeomorphism $h$ (the support map) between $\Omega'$ and $\Omega$ implies that $p=q$ (see for instance \cite[p. 120]{Fk}).

        Note that if $f \in
\noe$ and $y \in \Omega'$
satisfy \[{\partial^{ \lambda } f}
(h(y)) =0\] for all $\lambda \in \Lambda$, then by Lemma~\ref{46} we
have that \[(Tf)  (y)  =0.\]

        Now take $y \in \Omega'$ and fix ${\bf e} \in E$, ${\bf e}
\neq 0$. If $\# \Lambda$ stands for the
cardinal of $\Lambda$, then we can define a linear map $S_y
: {\Bbb R}^{\# \Lambda} \ra F$ as follows. Given
$(a_{\lambda_1}, 
\ldots , a_{\lambda_{\# \Lambda}} ) \in {\Bbb R}^{\#
\Lambda}$, we consider any $f \in \rno$ such that \[
{\partial^{\lambda } f} (h(y)) = a_{\lambda}\] for every $\lambda \in
\Lambda$. Then we define \[S_y (a_{\lambda_1},
\ldots , a_{\lambda_{\# \Lambda}} ) := (Tf {\bf e}) (y).\] The map $S_y$
is linear and, as we have seen above,
does not depend on
the function $f$ we choose. This implies that it is well defined.

        Then it is easy to see that there exist functions
$\alpha_{\lambda} $ from $\Omega'$ into
$F$, $\lambda \in \Lambda$,
such that for every $y \in \Omega'$ and every $f \in \rno$, 
\begin{equation}\label{for}
(Tf {\bf e}) (y) = \sum_{\lambda \in \Lambda} \alpha_{\lambda} (y)
{\partial^{ \lambda }f} (h(y)).
\end{equation}

        From now on we consider $i \in \{ 1, 2 , \ldots, p\}$
fixed.  Also $x_i \in \rno$ will be the projection
on the $i$-th coordinate for every point in $\Omega$.

        Next we define some
functions \[\alpha_i^{0}, \alpha_i^1 ,
\ldots , \alpha_i^{n+1}\] from $ \Omega' \times {\Bbb R}$ into $F$.
For every $y \in
\Omega'$ and $ t \in {\Bbb R}$,
\begin{eqnarray*}
\alpha_i^0 (y,t) &:=& (T\widehat{{\bf e}}) (y),\\ \alpha_i^1
(y,t) &:=& (Tx_i {\bf e}) (y) - \alpha_i^0 (y,t) t , \\ 2 ! \hs
\alpha_i^2 (y,t) 
&:=& (Tx_i^2 {\bf e}) (y) - \alpha_i^0 (y,t) t^2 -2 \alpha_i^1
(y,t) t ,\\ 3! \hs 
\alpha_i^3 (y,t) &:=& (Tx_i^3 {\bf e}) (y) - \alpha_i^0 (y,t)
t^3 -3 \alpha_i^1 
(y,t) t^2 - 6 \alpha_i^2 (y,t) t ,
\end{eqnarray*}
and in general, for $k \in \{1, 2, \ldots, n, n+1\}$
\begin{eqnarray*}
k! \hs
\alpha_i^k (y,t) &:=& (Tx_i^k {\bf e})
(y) - \alpha_i^0 (y,t) t^k - k \alpha_i^1 (y,t) t^{k-1} -k (k-1) \alpha_i^2
(y,t) t^{k-2} \\ &-& \ldots - k! \hs \alpha_i^{k-1} (y,t) t .
\end{eqnarray*}

\medskip

        {\bf Claim 1.} {\em  For $l \in \{0, 1, \ldots , 
n+1\}$, $l! \hs
\alpha_i^l $ is a polynomial in $t$ whose
coefficients are a linear combination of $T \widehat{{\bf e}},
Tx_i {\bf e} , \ldots,
Tx_i^l {\bf e}$. Moreover, for $y \in \Omega'$ fixed, the degree of the polynomial $l! \alpha_i^l (y, t)$ is at most $l$. 
If we also assume that $(T\widehat{{\bf e}}) (y)
\neq 0$, then the degree of $l! \alpha_i^l (y, t)$  is $l$
 and its leading
coefficient is equal to $ (-1)^l \alpha_i^0 (y,t)$ (notice that
this term does not
depend on $t$).}

        We are going to prove it by applying induction on $l$. It is clear that
this is true for
$l=0$. Suppose that this relation also holds for $l \in \{0, 1, 2,
\ldots, k\}$ for some $k \le n$. We are going to see that it holds for
$l=k+1$. We have that
\begin{eqnarray*}
 (k+1)! \hs
\alpha_i^{k+1} (y,t) &:=& (Tx_i^{k+1} {\bf e}) (y) - \alpha_i^0 (y,t) t^{k+1}
\\&-&
(k+1) \alpha_i^1 (y,t) t^k - (k+1) k \alpha_i^2 (y,t) t^{k-1} \\ &-&
(k+1) k (k-1) \alpha_i^3 (y,t) t^{k-2} \\ &-& \ldots \\ &-&
(k+1)! \hs \alpha_i^k
(y,t) t ,
\end{eqnarray*}
which implies that it is a
polynomial in $t$ and, for fixed $y \in \Omega$,  its
coefficient for the term $t^{k+1}$  is
$\alpha_i^0 (y,t) (-1 + (k+1) -
{k+1 \choose 2} + {k+1 \choose 3} - \ldots - (k+1) (-1)^{k} )$, which
is equal to $\alpha_i^0 (y,t) (-(1-1)^{k+1} + {k+1 \choose k+1}
(-1)^{k+1})$, that
is, to $(-1)^{k+1} \alpha_i^0 (y,t)$.

        Thus the claim is proved.

\medskip

Now, by Lemma~\ref{bueno},
we have that for every $k \in \{0, 1, \ldots, n +1\}$, $\alpha_i^k$
belongs to $C^m (\Omega' \times {\Bbb R}, F)$.

Next define $\alpha_0 :=T\widehat{{\bf e}}$, and for $k \in \{1, 2,
\ldots, n\}$, $\alpha_k := 
\alpha_{\lambda} : \Omega' \ra F$, where \[\lambda =
(\underbrace{0, 0, \ldots,0 , k,}_{i} 0, 
\ldots, 0).\] Also, let $h_i$ stand for the $i$-th coordinate
function of $h$.

\medskip

        {\bf Claim 2.}  {\em For every $k \in \{0, 1, \ldots ,
n\}$,
and for every $y \in \Omega'$, \[\alpha_k (y) = \alpha_i^k (y, h_i (y)).\]}

First we have from Equation~\ref{for} that for $k
\in \{0, 1, \ldots, n\}$ and $y \in \Omega'$,
\[(Tx_i^k {\bf e}) (y) = \sum_{\lambda \in \Lambda} \alpha_{\lambda} (y)
{\partial^{\lambda } x_i^k} (h(y)),\]
which can be written as

\begin{eqnarray}\label{dafa}
(Tx_i^k {\bf e}) (y) &=& \alpha_0 (y) h_i^k (y) \\&+&
 k \alpha_1 (y)
h_i^{k-1} (y) \nonumber \\  &+&  \ldots \nonumber \\&+& k! \alpha_{k-1} (y)
h_i(y) \nonumber \\  &+&  k!
\alpha_k (y), \nonumber
\end{eqnarray}
because \[{\partial^{ \lambda }
x_i^k}
(h(y)) =0\] whenever $\lambda \in \Lambda$, $\lambda \neq (\underbrace{0, 0,
\ldots , j ,}_{i} 0, \ldots, 0)$, $j \in \{0, 1, \ldots , k\}$.

        On the other hand, it is clear that $\alpha_0 (y) = (T
\widehat{{\bf e}}) (y) = \alpha_i^0 (y, h_i (y))$ for every $y \in \Omega'$. Also suppose that $k <n$ and that $\alpha_j (y) = \alpha_i^j (y, h_i
(y))$ for every $j \in \{0,1, \ldots ,k\}$ and every $y \in
\Omega'$. Then, by Equation~\ref{dafa}, for $y \in \Omega'$
\begin{eqnarray*}
(k+1)! \alpha_{k+1} (y) &=& (Tx_i^{k+1} {\bf e}) (y) \\&-&
\alpha_0 (y) h_i^{k+1}  (y) \\&-& (k+1) \alpha_1 (y)
h_i^{k} (y)\\ &-& \ldots \\&-& (k+1)! \alpha_{k} (y) h_i(y),
\end{eqnarray*}
which coincides with
$(k+1)! \alpha_i^{k+1} (y, h_i (y))$, and the claim is proved.

\medskip

On the other hand, notice that in the same way as we obtain Equation~\ref{dafa}, we have

\begin{eqnarray}\label{das}
(Tx_i^{n+1} {\bf e}) (y) &=& \alpha_0 (y) h_i^{n+1} (y) \\&+&
(n+1) \alpha_1 (y) 
h_i^{n} (y) \nonumber \\&+& \ldots \nonumber \\&+& (n+1)!
\alpha_{n} (y) h_i(y) \nonumber 
\end{eqnarray}

for every $y \in \Omega'$.

\medskip

        {\bf Claim 3.}  {\em Suppose that $y_0 \in \Omega'$ satisfies
$\alpha_0 (y_0) \neq 0$. Then for every open 
neighborhood $U$ of $y_0$, there
exists a nonempty open subset $W'_i$ of $U$ where $h_i$
is of class $C^m$.}

        First we define
$F_i^0: \Omega' \times {\Bbb R} \ra F$
as $F_i^0:=  (n+1)! \alpha_i^{n+1}$, that is,
\begin{eqnarray*}
F_i^0(y,t) &=& (Tx_i^{n+1} {\bf e}) (y) \\&-& \alpha_i^0 (y,t)
t^{n+1} \\&-&(n+1) 
\alpha_i^1 (y,t) t^n \\ &-& (n+1) n \alpha_i^2 (y,t) t^{n-1} \\&-&
\ldots \\&-& (n+1)! \hs \alpha_i^n (y,t) t,
\end{eqnarray*}
for every $y \in \Omega'$, $t \in {\Bbb R}$.

        Then, if $j \in \{1, 2, \ldots , n+1\}$, we define
$F_i^j: \Omega' \times
{\Bbb R} \ra F$ as \[F_i^j (y,t) = \frac{\partial^j F_i^0}{\partial t^j}
 (y,t),\] for all $y \in \Omega'$, $t \in {\Bbb R}$.

Notice that from the definition of $F_i^0$, Claim 2 and
Equation~\ref{das}, we deduce that \[F_i^0 (y, h_i (y)) =0\] for 
every $y \in \Omega'$. 
Also, as we stated in Claim 1, the coefficients of $F_i^0$ as a
polynomial of degree $n+1$ in $t$ are
linear combinations of \[(T \widehat{{\bf e}})(y), (Tx_i
{\bf e}) (y), \ldots, 
(Tx_i^{n+1}  {\bf e}) (y),\] and consequently, by
Lemma~\ref{bueno}, for $k \in \{1, 2, \ldots ,
n\}$, $F_i^k$ belongs to
$C^m (\Omega' \times {\Bbb R}, F)$.

Taking into account that $F_i^{n+1} (y, t) = (n+1)! (-1)^{n+1} \alpha_0 (y)$ for
every $(y, t) \in \Omega' \times {\Bbb R}$, and the fact that 
$\alpha_0 (y_0) \neq 0$, there
exists $k_0 \in \{0, 1, \ldots, 
n\}$ such that $ F_i^{k_0} (y ,h_i (y ))  
= 0$ for every $y$ in a neighborhood of $y_0$ and $ F_i^{k_0 +1}
(y ,h_i (y ))   $ takes a value different from $0$ for some $y$
in every neighborhood of $y_0$. Suppose then that $U$ is an open
neighborhood of $y_0$ such that $ F_i^{k_0} (y ,h_i (y ))  
= 0$ for every $y \in U$ and that $y_1 \in U$ satisfies $
F_i^{k_0 +1 } (y_1 ,h_i (y_1 ))  = {\bf f} \in F$, ${\bf f}
\neq 0$. Now take ${\bf f}'$ in the dual space $F'$ (where $F$
is wiewed as a {\em real} Banach space) such that
${\bf f}' ({\bf f}) \neq 0$. According to the Implicit
Function Theorem (\cite[p.148]{Fl}), there 
exist a neighborhood $V$ of $(y_1 , h_i (y_1))$, an open
neighborhood $W$ of $y_1$, and a function $\phi : W \ra {\Bbb R}$
of class $C^m$ such that $\phi (y_1) = h_i (y_1)$ and \[\{ (y, t) \in
V : {\bf f}' \circ F_i^{k_0} (y,t) =0 \} = \{ (y, \phi (y) ) :
y \in W \}.\]
It is easy to prove that this implies that $\phi \equiv h_i$ on
a neighborhood $W'_i$ of $y_1$, that is, for every open
neighborhood $U$ of $y_0$, there
exists an open subset $W'_i$ of $U$ where $h_i$
is of class $C^m$.
The claim is proved.

\medskip

Since both $T$ and $T^{-1}$ are biseparating, we can assume from now on, without loss of generality, that $n \le m$.

\medskip

{\bf Claim 4.} {\em Suppose that $U$ is a nonempty open subset of
$\Omega'$. Then there exists a nonempty open subset $W'$ of $U$ such that
the restriction of $h$ to $W'$ is a diffeomorphism of class $C^m$.}

Notice first that the the open set $\{y \in \Omega': \alpha_0 (y) \neq 0\}$ is
dense in $\Omega'$. Otherwise we could find $g \in \moe$, $g \neq 0$, such that
$ c(T \widehat{{\bf e}}) \cap  c(g) =  c(\alpha_0) \cap  c(g) =
\emptyset$. Since $T^{-1}$ is  separating, this would give us $\Omega \cap c(T^{-1}g) =
c(\widehat{{\bf e}}) \cap c(T^{-1} g ) = \emptyset$, which is
impossible.

        So far we have considered $i \in \{1, 2, \ldots, p\}$
fixed. Of course a similar process can be done for every $i \in
\{1, 2, \ldots, p\}$. In particular, taking into account the above
paragraph, by Claim 3, there exists an
open 
subset $W'_1$ of $U$ such that $h_1$ is of class $C^m$ in
$W'_1$. For the same reason we can find an open subset $W'_2$ of
$W'_1$ where $h_2$ is of class $C^m$. Following this process we
construct (nonempty) open sets $W'_1, \ldots, W'_p$ with $W'_1
\supset W'_2 \supset \ldots \supset W'_p$ such that $h$ is
of class $C^m$ in $W'_p$. It is clear that a similar reasoning
shows that the map 
$h^{-1}$ is of class $C^n$ in an open subset $V$ of $h(W'_p)$.
Then our situation is as follows: $h^{-1}$ is of class $C^n$ in
$V$ and $h$ is of class $C^m$ in $h^{-1} (V)$. It is well known
that, since $m \ge n \ge 1$, this implies that $h$ is a
diffeomorphism of class $C^m$ in 
$W':= h^{-1} (V)$, and we are done.

\medskip

        {\bf Claim 5.} {\em Let $W'$ be as in Claim 4. For every $k
\in \{0, 1, \ldots, n\}$, 
the map $\alpha_k$ belongs to
$C^m (W', F)$.}

 First we have that $\alpha_0 = T \widehat{{\bf e}}$ belongs to $\moe$. It is
also clear that if $k \in \{0, 1, \ldots, 
n-1\}$, then as given in Equation~\ref{dafa}, on $\Omega'$,
\[ (k+1) ! \hs \alpha_{k+1} =
T x_i^{k+1} {\bf e} - \alpha_0 h_i^{k+1} - 
(k+1) \alpha_1 h_i^{k} - \ldots - (k+1) ! \hs \alpha_k
h_i .\] 

Consequently, since $h_i \in C^m (W')$, if $\alpha_0, \alpha_1, 
\ldots, \alpha_k$
belong to $C^m (W', F)$, $\alpha_{k+1}$ also belongs to $C^m
(W', F)$ and then we are done.

\medskip

        {\bf Claim 6.} {\em For every $k \in \{1, 2, \ldots, n\}$,  $\alpha_{k} \equiv 0$ in $W'$.}

       Suppose that $y_0 \in W' $ and $\alpha_n (y_0)
\neq 0$. Since by Claim 5 $\alpha_n$ is continuous in $W'$, there exists an open   neighborhood $U(y_0)$ of $y_0$ such that $U(y_0) \subset W'$ and $\alpha_n (y) \neq 0$ for every $y \in U(y_0)$. Then take $g \in C^n ({\Bbb R}, {\Bbb R})$ such
that $g^{(n)}$ is not 
derivable at the point $h_i (y_0)$. We define $f \in \rno$
as \[f (x) := g(x_i)\] for every $x = (x_1 , x_2 , \ldots , x_p)
\in \Omega$. In this way we have that
\[\frac{\partial^{n+1}f}{\partial x_i^{n+1}} (h(y_0))\]
does not exist. Consequently, using a reasoning similar to that
giving Equation~\ref{dafa}, we have that
Equation~\ref{for} applied to 
$f$ is
\[(Tf{\bf e}) (y) = \alpha_0 (y) f(h(y)) + \alpha_1 (y)
\frac{\partial f}{\partial x_i} 
(h(y)) + \ldots + \alpha_n (y) \frac{\partial^n f}{\partial x_i^n}
(h(y))\]for every $y \in U(y_0)$. 
Now we analyse the terms in the above equation, taking into account that we are assuming  $1 \le n \le m$, and that by Claim 4, $h$ is a diffeomorphism of class $C^m$ in $W'$. First $Tf{\bf e}$ is of class $C^1$ in $U(y_0)$. Also, by Claim 5, for $k \in \{0, 1, \ldots ,n\}$, each $\alpha_k$ is of class $C^1$ in $U(y_0)$. Finally,  $f$ and all of its partial derivatives up to  order $n-1$ are of class $C^1$ in $\Omega$.

Thus we deduce from the above equation that
\[ \alpha_n \frac{\partial^{n}  f}{\partial x_i^{n}}\circ h\]
is of class $C^1$  in $U(y_0)$, and clearly the same applies to the function \[ \frac{\partial^{n}  f}{\partial x_i^{n}}\circ h,\] by Lemma~\ref{dabute}.

But, as we said before, $h$ is a diffeomorphism of class $C^m$ in $W'$, and
consequently \[\frac{\partial^{n}  f}{\partial x_i^{n}}\]
admits a partial 
derivative with respect to the $i$-th coordinate at the point $h(y_0)$,
which is a contradiction. This implies that $\alpha_n \equiv 0$
in $W'$. In a
similar way we can see that $\alpha_k \equiv 0$ in $W'$ for $k
\ge 1$, that is, $\alpha_{\lambda} \equiv 0$ in $W'$ whenever
$\lambda \in \Lambda$ is of the form
\[\lambda = (\underbrace{0, 0, \ldots , k,}_{i} 0,
\ldots, 0).\]

\medskip

        {\bf Claim 7.} {\em For every
$\lambda \in \Lambda - \{(0, 0, 
\ldots, 0)\}$, $\alpha_{\lambda} \equiv 0$ in $W'$.}

Here our reasoning will be similar to the one given in Claim 6.  In this way, if $i , j \in \{1, 2, \ldots, p\}$, $i 
\neq j$, 
again by Equation~\ref{for}, for every $y \in W'$, \[ (T (x_i x_j {\bf e}) )(y) =
\alpha_0 (y) h_i (y) h_j (y) + \alpha_{\lambda_0^1} (y)\] where
$\lambda_0^1 := (\lambda_{01}^1 , \lambda_{02}^1, \ldots,
\lambda_{0p}^1 )$, $\lambda_{0i}^1 =1 =\lambda_{0j}^1$ and $\lambda_{0k}^1
=0$, whenever $k \neq i, j$. Taking into account that
$\alpha_0$, $h_i$ and $h_j$ are of class $C^m$ in $W'$, we easily
deduce that $ \alpha_{\lambda_0^1} $ is of class $C^m$ in $W'$.
Likewise, we can inductively prove that
$\alpha_{\lambda_0^l}$ is of class $C^m$ in $W'$, where
$\lambda_0^l := (\lambda_{01}^l , \lambda_{02}^l, \ldots,
\lambda_{0p}^l )$, $\lambda_{0i}^l =l$, $\lambda_{0j}^l =1$ and $\lambda_{0k}^l
=0$, whenever $k \neq i, j$, $l \in \{1, 2, \ldots, n-1 \}$.
Suppose that
$\alpha_{\lambda_0^{n-1}} (y_0) \neq 0$ for some $y_0 \in W'$. Then, as in the proof of Claim 6, we take an open  neighborhood $U(y_0)$ of $y_0$ such that $U(y_0) \subset W'$ and $\alpha_{\lambda_0^{n-1}}  (y) \neq 0$ for every $y \in U(y_0)$.

Also, we take 
$f (x) = g(x_i)$ for every $x = (x_1 , x_2 , \ldots , x_p)
\in \Omega$, where  these functions meet the same requirements as in the proof of Claim 6,
and define \[d(x) := x_j f(x),\]
for every $x \in \Omega$. Clearly $d$ just depends on
the $i$-th and $j$-th
coordinates, which implies that its only partial derivatives
which possibly are not zero at $h(y) \in h(U(y_0))$ are maybe those
\[{\partial^{ \lambda } d} \]
for
\[\lambda = (\underbrace{0, 0, \ldots , 1,}_{j} 0,
\ldots, 0),\]
\[\lambda = (\underbrace{0, 0, \ldots , k,}_{i} 0,
\ldots, 0),\] $k=1, 2, \ldots, n$, or
\[\lambda = \lambda_0^l,\] $l \in \{1,2,
\ldots, n-1\}$. Taking into account that $\alpha_{\lambda}
\equiv 0$ on $W'$ for 
\[\lambda = (\underbrace{0, 0, \ldots , 1,}_{j} 0,
\ldots, 0)\] and
\[\lambda = (\underbrace{0, 0, \ldots , k,}_{i} 0,
\ldots, 0),\] Equation~\ref{for} gives us, for every $y \in U(y_0)$,
\[ (T d {\bf e} )(y) =
\alpha_0 (y) d(h(y)) + \alpha_{\lambda_0^1 } (y) \frac{\partial
f}{\partial x_i} ( h(y)) + \ldots + \alpha_{\lambda_0^{n-1}} (y)
\frac{\partial^{n-1} f}{\partial x_i^{n-1}}
(h(y)).\] We deduce as in the proof of Claim 6 that
\[ \frac{\partial^{n-1}
f}{\partial x_i^{n-1}} ( h(y)) \] admits a second partial derivative
with respect to $x_i$ at the point $h(y_0)$, which is a contradiction. This implies
that $\alpha_{\lambda_0^{n-1}} \equiv 0$ in $W'$. In the same
way we deduce that $\alpha_{\lambda_0^{l}} \equiv 0$ in $W'$,
for $l \in \{1, 2, \ldots, n-2\}$.

A similar pattern of proof leads
us to the fact that $\alpha_{\lambda} \equiv 0$ in $W'$ for every
$\lambda \neq (0, 0, \ldots ,0)$, $\lambda \in \Lambda$.

\medskip

{\bf Claim 8.}  {\em For every
$\lambda \in \Lambda - \{(0, 0, 
\ldots, 0)\}$, $\alpha_{\lambda} \equiv 0$ in $\Omega'$.}

        Notice that from an open subset $U$ of $\Omega'$, in Claim
4 we obtain a subset $W'$ of $U$. Notice also that this process can be
done for any open subset of $\Omega'$ because, as we saw in the
proof of Claim 4, $c(\alpha_0)$ is dense in $\Omega'$. Also in
Claim 7 we proved 
that, for $\lambda \neq (0, 0, \ldots, 0)$,  $\alpha_{\lambda}
\equiv 0$ on all the subsets $W'$ obtained in this way. This
implies clearly that all these functions $\alpha_{\lambda}$ are
equal to $0$ on a dense subset of $\Omega'$. Consequently, to
prove Claim 8, it is enough to show that all these functions are
continuous.

        We are going to prove it using induction on
$\va \lambda \vb$. First, for $\va \lambda \vb =0$, we have that $
\alpha_{(0, \ldots, 0) } = T \widehat{{\bf e}} $ belongs to
$\moe$ and, consequently, it is continuous.

        Now assume that $k \le n-1$, and whenever $\va \lambda \vb
\le k$, then $\alpha_{\lambda}$ is a continuous function.
Then fix $\lambda = (\lambda_1 , \ldots , \lambda_p ) \in \Lambda$ with
$\va \lambda \vb = k+1$.

        Next define $f \in \rno $ as \[f := x_1^{\lambda_1} 
x_2^{\lambda_2}  \ldots x_p^{\lambda_p} .\]

        It is clear that given $\mu = (\mu_1 , \ldots , \mu_p)
\in \Lambda$, if $\mu_i > \lambda_i$ for some $i \in \{1, \ldots
, p\}$, then \[{\partial^{ \mu
} f}  (h(y)) =0\] for every $y \in \Omega'$. This implies that in our situation, 
Equation~\ref{for} can
be written as 
\[(Tf {\bf e}) (y) = \sum_{\mu << \lambda} \alpha_{\mu} (y)
{\partial^{ \mu }f}  (h(y)),\]
 where $\mu << \lambda$
means $\mu_i \le \lambda_i$ for every $i \in \{1, \ldots , p\}$.

        As a consequence, for every $y \in \Omega'$, 
\[ \lambda_1 ! \ldots \lambda_p ! \alpha_{\lambda} (y) =
\alpha_{\lambda} (y) {\partial^{\lambda }f} (h(y)) = (Tf {\bf e}) (y) - \sum_{\mu << \lambda,
\mu \neq \lambda} \alpha_{\mu} (y) 
{\partial^{ \mu }f}  (h(y)).\]

        On the other hand, if $\mu << \lambda$ and $\mu \neq
\lambda$, then $\va \mu \vb \le k$, and consequently, taking
into account that $h$ is continuous and the hypothesis of
induction, we deduce that $\alpha_{\lambda}$ is continuous, and the claim is proved.

\medskip

        Recall that all the process developed so far concerns
functions of the form $f {\bf e} \in \noe$, where ${\bf e} \in E-\{0\}$
and $f \in \rno$. For this ${\bf e}$, we define \[a_{{\bf e}}
:= \alpha_0 = T \widehat{{\bf e}}.\] Notice that by Claim 8, we have
\begin{equation}\label{asube}
(Tf {\bf e}) (y) = a_{{\bf e}} (y) f (h(y))
\end{equation}
 for
every $y \in \Omega'$ and every $f \in \rno$.

        Next we define a map $J: \Omega' \ra L'(E,F)$ as $(Jy)
(0) =0$, and \[(Jy)
({\bf e}) := a_{{\bf e}} (y)\] for each $y \in \Omega'$ and ${\bf e}
\in E- \{0\}$.

\medskip

        {\bf Claim 9.} {\em For every $f \in \noe$ and $y \in
\Omega'$, \[(Tf) (y) = (Jy) (f(h(y))).\]}

        Fix $y \in \Omega'$. Suppose that \[
{\partial^{ \lambda } f} (h(y)) = {\bf e}_{\lambda} \in E\] for each
$\lambda \in \Lambda$.

       Let $\Lambda_* := \{\lambda \in \Lambda : {\bf e}_{\lambda} \neq 0 \}$. Next,  for each $\lambda \in \Lambda_*$, take a
function $f_{\lambda} \in \rno$ such that \[
{\partial^{\lambda } f_{\lambda}} (h(y)) =1\] and \[
{\partial^{ \mu } f_{\lambda}} (h(y)) = 0\] for every $\mu \neq \lambda$, $\mu \in \Lambda$.

        It is easy to see that, for every $\mu \in \Lambda$,
\[
{\partial^{ \mu} \sum_{\lambda \in \Lambda_*}
f_\lambda {\bf e}_{\lambda}} (h(y)) = {\bf e}_{\mu},\]
if $\mu \in \Lambda_*$, and
\[
{\partial^{ \mu} \sum_{\lambda \in \Lambda_*}
f_\lambda {\bf e}_{\lambda}} (h(y)) = 0\]
if $\mu \notin \Lambda_*$.

        According to Lemma~\ref{46}, this implies that \[(Tf)
(y) = \left( T \sum_{\lambda \in \Lambda_*} f_{\lambda}
{\bf e}_{\lambda} \right) ( y ).\] Consequently, by Equation~\ref{asube},
\[(Tf) (y) = \sum_{\lambda \in \Lambda_*} a_{{\bf e}_{\lambda}}
(y) f_{\lambda} (h(y)).\]

        But by the way we have constructed the functions
$f_{\lambda}$, we have that $f_{\lambda} (h(y)) =0$ if $\lambda
\neq (0,0, \ldots, 0)$. 

On the other hand, let us denote by
$0$ the multiindex $(0, 0, \ldots, 0)$. If $0 \notin \Lambda_*$, that is, if ${\bf e}_0 =0$, we conclude from the above equality that $(Tf)(y)=0= a_{{\bf e}_0} (y) = (Jy) (0)$. Finally, if $0 \in \Lambda_*$, taking into account that $f_0 (h(y)) =1$,  we deduce that
\[(Tf) (y) = a_{{\bf e}_0} (y) f_0 (h(y)) = a_{{\bf e}_0} (y)
= (Jy) ({\bf e}_0),\]
and we are done.

\medskip

        {\bf Claim 10.} {\em Given $y \in \Omega'$, there exists
${\bf e} \in E$ such that $a_{{\bf e}} (y) \neq 0$.}

        Notice that $T$ is bijective, so if ${\bf f}
\in F$, ${\bf f} \neq 0$, there exists $g \in
\noe$ with $Tg = \widehat{{\bf f}}$. In particular, by Claim 9,
we have that $(Jy) (g(h(y))) = {\bf f}$. In other words, if we
take $y \in c(g \circ h)$ and define
${\bf e} := g(h(y))$, we have that $a_{{\bf e}} (y) = (Jy)({\bf e})=
{\bf f} \neq 0$.

\medskip

        {\bf Claim 11.} {\em $h$ is a function of class $C^m$.}

        Fix $y_0 \in \Omega'$. By Claim 10, we can take
${\bf e} \in E$ such that $a_{{\bf e}} (y)  \neq 0$. Since
$a_{{\bf e}} = T \widehat{{\bf e}}$, it is a 
continuous function, and we deduce that for every $y$ in some neighborhood
$V$ of $y_0$, $a_{{\bf e}} (y) \neq 0$.

        Now recall that Equation \ref{asube}, 
\[(Tf {\bf e}) (y) = a_{{\bf e}} (y) f (h(y)),\]
holds in particular for
every $y \in V$ and every $f \in \rno$.

Consequently, for  $i \in \{1, 2, \ldots ,
p\}$, 
\[(Tx_i {\bf e}) (y) = a_{{\bf e}} (y) h_i (y)\] for every $y
\in V$. Since $a_{{\bf e}}
(y) \neq 0$ for every $y \in
V$, applying Lemma~\ref{dabute}, we have that $h_i$ is of
class $C^m$ in $V$. Clearly this implies that $h$ is  of class
$C^m$, and we are done.

\medskip

        {\bf Claim 12.} {\em $n=m$ and $h$ is a diffeomorphism of class $C^m$.}

        Recall that we are assuming that $n \le
m$. Now, we have that by Claim 4, for every nonempty open set $V \subset
\Omega'$, there is a nonempty open set $V' \subset V$ such that
the restriction of $h$ to $V'$ is a diffeomorphism of class
$C^m$. 
Take $y_0 \in \Omega'$. By Claim 10, there exists ${\bf e} \in
E$ and an open  neighborhood $V$ of $y_0$ such that $a_{{\bf e}}
(y) \neq 0$ for every $y \in V$. Now, as we mentioned above,
there exists an open set $V'$, $V' \subset V$, where the restriction of
$h$ is a diffeomorphism of class $C^m$. Assume now that $n <m$ and take $g \in
A^n (\Omega', {\Bbb R} ) - C^m (\Omega', {\Bbb R})$ such
that $c(g) \subset V'$. Next
define $f : \Omega \ra E$, as \[f(x) := g 
(h^{-1} (x)) {\bf e},\] for each $x \in \Omega$. We are going
to prove that $g$ is of 
class $C^m$, obtaining a contradiction.

        It is immediate that $f \in \noe$, and applying
Equation~\ref{asube}, we get
\[(Tf) (y) = a_{{\bf e}} (y) g(h^{-1} (h(y))),\] that is, 
\[(Tf) (y) = a_{{\bf e}} (y) g (y),\] for every $y \in \Omega'$.

Now we have that $a_{{\bf e}} (y) \neq 0$ for every $y
\in V$. Finally, by Lemma~\ref{dabute}, $g$ is of class
$C^m$ in $V$, and so is in $\Omega'$, which contradicts our
assumption.
This implies that $m \le n$, and since we are assuming that
$n \le m$, the claim
is proved.

\medskip

        {\bf Claim 13.} {\em For every $y \in \Omega'$, $Jy \in L'(E,F)$
is bijective.}

        Since $m=n$, all claims above also hold for
$T^{-1}$, and this means that there exists $K: \Omega \ra
L'(F,E)$ such that for every $g \in \moe$ and $x \in \Omega$,
$(T^{-1} g) (x) = (Kx) (g(h^{-1} (x)))$. 

        Fix $y \in \Omega'$ and ${\bf f} \in F-\{0\}$. Let $x =
h(y)$. Now take $g \in 
\moe$ with $g(y) = {\bf f}$. Then its is clear that ${\bf f}
=g (y) =  (T (T^{-1} g)) (y)$, that is,
\begin{eqnarray*}\label{reves}
{\bf f} &=& (Jy) ((T^{-1} g) (x))\\
&=& (Jy) ((Kx) (g(h^{-1} (x))))\\
&=& (Jy) ((Kx) (g(y)))\\
&=& (Jy) ((Kx) ({\bf f})).
\end{eqnarray*}

        This implies that $(Jy) (Kx)$ is the identity map on
$F$. In the same way we can prove that $(Kx)(Jy)$ is the
identity map on $E$. Consequently, $Jy$ is bijective.

\medskip

This ends the proof of the lemma.
\end{pf}

\section{Some results on automatic continuity}

In this section we follow the same notation as in the previous one.

Here, if we endow the spaces with some natural topologies, then we obtain the continuity as a consequence. Notice that, according to Lemma~\ref{401}, we can assume in particular that $n=m$ and $p=q$.

\begin{definition}\label{addi}
We say that a locally convex topology in $\noe$ is {\em compatible with the pointwise convergence} if the following two conditions are satisfied:
\begin{itemize}
\item 1. when endowed with it, $\noe$ is a Fr\'echet (or Banach) space, and
\item 2. if $(f_n)$ is a  sequence in $\noe$ converging to zero, then $(f_n (x))$ converges to zero for every $x \in \Omega$.
\end{itemize}
\end{definition}

\begin{theorem}\label{401b}
Assume that $\noe$ and $\nmoe$ are endowed with any topologies which are compatible with the pointwise convergence. Suppose that $T: \noe \ra \nmoe$
is a ${\Bbb K}$-linear  biseparating map. Then $T$ is continuous.
\end{theorem}

\begin{pf}
In our proof we will take advantage of the description of $T$ given in Lemma~\ref{401}. For this reason we will use the notation given there. We start proving the following claim.

\medskip

        {\bf Claim.} {\em Let $U$ be a (nonempty) bounded open subset of $\Omega'$ with $\clp U \subset \Omega'$. Then the set \[A:= \{y \in U : Jy \in B'(E, F) \mathrm{\hspace{.03in} is \hspace{.03in} not \hspace{.03in} continuous}\}\] is finite.}

	Suppose that this is not the case, but there exist infinitely many $y \in U$ such that $Jy$ is not continuous. 
We are going to construct inductively a sequence of points in $A$, a sequence $(U_n)$ of pairwise disjoint open subsets of $U$, a sequence of functions $(f_n)$ in $\noor$, and a sequence of norm-one elements of $E$, satisfying the following properties:

\begin{itemize}
\item 1. $ h(y_n) \in c(f_n) \subset h(U_n)$ for every $n \in {\Bbb N}$.
\item 2. $\vc f_n \vd: = \max_{\lambda \in \Lambda} \sup_{x \in \Omega} \va \partial^{\lambda} f_n (x) \vb =1/2^n$ for every $n \in {\Bbb N}$.
\item 3. $\vc (Jy_n) ({\bf e}_n) \vd \ge n/\va f_n (h(y_n))\vb$ for every $n \in {\Bbb N}$.
\end{itemize}

Take any point $y_1 \in A$ such that there are accumulation points of $A$ in $A - \{y_1\}$. Then consider an open subset $U_1$ of $U$ in such a way that $y_1 \in U_1$, and there are infinitely many points of $A$ outside $\clp U_1$. Next take $f_1 \in \noor$ such that $\vc f_1 \vd =1$, and such that $h(y_1) \in c(f_1) \subset h(U_1)$. Since $Jy_1$ is not continuous, there exists ${\bf e} \in E$, $\vc {\bf e } \vd =1$, with \[\vc (J y_1)({\bf e}) \vd \ge \frac{1}{\va f_1 (h(y_1)) \vb}.\]

Now assume that we have $\{y_1 , y_2 , \ldots , y_n\} \subset A$, $U_1, U_2 , \ldots, U_n \subset U$ open and pairwise disjoint such that there are infinitely many points of $A$ outside $\clp U_1 \cup \clp U_2  \cup \ldots \cup \clp U_n$, $\{f_1 , f_2 , \ldots , f_n\} \subset \noor$ with $h(y_i) \in c(f_i) \subset h(U_i)$ and $\vc f_i \vd =1/2^i$, for $i \in \{1, 2, \ldots, n\}$, and ${\bf e}_1, {\bf e}_2 , \ldots, {\bf e}_n \in E$ all of them with norm $1$, and such that $\vc (Jy_i) ({\bf e}_i) \vd \ge i
/\va f_i (h(y_i))\vb$ for $i = 1, 2, \ldots, n$.

Now it is easy to see how to take $y_{n+1}$, $U_{n+1}$, $f_{n+1}$, and ${\bf e}_{n+1}$ so that Properties 1, 2 and 3 above hold.

Since $\vc f_n \vd = 1/2^n$ for every $n \in {\Bbb N}$, we deduce that the map \[g := \sum_{n=1}^{\infty} f_n {\bf e}_n\] belongs to $\noe$. Consequently, $Tg$ should belong to $\nmoe$. But we know by Claim 9 in the proof of Lemma~\ref{401} that $(Tg) (y_n) = (Jy_n) (g(h(y_n)))$ for every $n \in {\Bbb N}$. This implies, by Property 3 above, 

\begin{eqnarray*}
\vc (Tg) (y_n) \vd &=& \vc (Jy_n) (f_n (h(y_n)) {\bf e}_n) \vd \\
&=& \va f_n(h(y_n)) \vb \vc (Jy_n) ({\bf e}_n) \vd \\
&\ge& n.
\end{eqnarray*}

As a consequence $Tg$ is unbounded in $U$. Since this is not possible, we conclude that the claim is correct.

\medskip

Next, it is clear that to prove that $T$ is continuous it is enough to show that it is closed, because we are dealing with Fr\'echet spaces.
To prove it, let us consider a sequence $(f_n)$ in $\noe$ convergent to zero, and assume that $(Tf_n)$ converges to $g \in \nmoe$. We are going to prove that $g=0$.

Take a  bounded open subset $U$ of $\Omega'$ with $\clp U \subset \Omega'$. By the claim above, we have that the subset $A$ of points $y \in U$ such that $Jy$ is not continuous is finite. So,  if $y \in U-A$, $Jy$ belongs to $B(E,F)$. Consequently, since $f_n(h(y))$ goes to zero (because the topology in $\noe$ is compatible with the pointwise convergence), then we have that $(Tf_n) (y) = (Jy) (f_n (h (y)))$
also goes to zero, that is, $g(y)$ must be zero. But  taking into account that $U-A$ is dense in $U$, we deduce that $g \equiv 0$ on $U$. The conclusion follows now easily and $T$ is continuous.
\end{pf}

\section{Biseparating maps and functions of class ${\rm s}-C^{n}$}

As in the previous section, in this one we also follow the notations introduced in Section 3. Our aim will be to give a final description of biseparating maps between spaces of vector-valued differentiable functions taking into advantage that we know that they must be continuous when the spaces are endowed with some natural topologies. Of course these topologies will be compatible with the pointwise convergence. Namely, it is well known that by means of the
seminorms $p_K$ defined as \[p_K
(f):= \max_{\lambda
\in \Lambda} \max_{x \in K} \vc \partial^{\lambda}
f (x) \vd \] for $f \in C^n (\Omega, E)$, where
$K$ runs through the compact subsets of $\Omega$, $C^n (\Omega, E)$
becomes a locally convex space. In fact it is a Fr\'echet space. 
In the same way, in $C^n
(\bar{\Omega}, E)$ we can consider the norm $\vc \cdot \vd$ defined as \[\vc f\vd := \max_{\lambda
\in \Lambda} \sup_{x \in \Omega} \vc \partial^{\lambda}
f (x) \vd \] for $f \in C^n
(\bar{\Omega}, E)$. With this norm, our space $C^n
(\bar{\Omega}, E)$ is also complete. We assume that $\noe$ and $\nmoe$ are endowed with the above topologies. Remark also that, as it follows easily from the Closed Graph Theorem, the topologies compatible with the pointwise convergence in our spaces coincide with these topologies.

Suppose that $K: \Omega \ra L(E,F)$ is a continuous map, where $L(E,F)$ is endowed with the topology of the norm. For each ${\bf e} \in E$, we define $K_{\bf e} : \Omega \ra F$ as $K_{\bf e} (y) := (Ky) ({\bf e})$ for every $y \in \Omega$. We say that $K$ is of class ${\rm s}-C^1$ if, for every ${\bf e} \in E$, the map $K_{\bf e}$ admits all partial derivatives of order $1$ in $\Omega$, and for each $i = 1, \ldots, p$, the map
\[\frac{\partial_{\rm s}}{\partial x_i} K : \Omega \ra L(E,F),\] sending each $y \in \Omega$ and each ${\bf e} \in E$ into $\frac{\partial}{\partial x_i} K_{\bf e} (y)$, is continuous when 
considering in $L(E,F)$ the strong operator topology, that is, the coarsest topology such that the mapping $A \in L(E,F) \hookrightarrow A{\bf e} \in F$ is continuous for every ${\bf e} \in E$.

\begin{definition}
Let $n \ge 2$. A map $J : \Omega \ra L(E,F)$ is said to be of class ${\rm s}-C^{n}$ if the following three statements are satisfied:
\begin{itemize}
\item 1. $J$ is of class $C^{n-1}$ (considering $L(E,F)$ as a Banach space).
\item 2. All partial derivatives $K: \Omega \ra L(E,F)$ of order $n-1$ of $J$ are of class ${\rm s}-C^{1}$.
\end{itemize}
\end{definition}

Next proposition states that when $J: \Omega' \ra L(E,F)$ is of class ${\rm s}-C^n$, we can define maps through $J$ from $A^n (\Omega, E) \ra A^n (\Omega', F)$ in a natural way.

\begin{prop}\label{dadada}
Suppose that $n=m$, $p=q$. Let $J: \Omega' \ra L(E,F)$ be a map of class  ${\rm s}-C^{n}$, and let $h$ be a diffeomorfism of class $C^n$ from $\Omega'$ onto $\Omega$.
If, for $f \in C^n (\Omega, E)$, we define  $(Tf)(y) := (Jy) (f(h(y)))$ for every $y \in \Omega'$, then $Tf \in C^n (\Omega' , F)$.
\end{prop}

\begin{pf}
We consider first the map
$\Phi : L(E,F) \times E \ra F$ defined as $\Phi (A , {\bf e}) :=
A {\bf e}$ for each $(A, {\bf e}) \in L(E,F) \times
E$. This is clearly bilinear and continuous  when we consider in $L(E,F)$ the topology of the norm.

Suppose next that $L(E,F)$ is endowed again with the topology of the norm, and that $K: \Omega' \ra L(E,F)$ is a continuous map. Then, given  $g: \Omega' \ra E$ continuous, the map $S_{K}^{g}: \Omega'
\ra L(E,F) \times E$ sending each $y \in \Omega'$ into $(Ky,
g(y))$ is continuous. 

On the other hand, if we suppose  that $L(E,F)$ is endowed with the topology of the norm and that $f \in C^n (\Omega, E)$, then $S_{J}^{f \circ h}$ is of class $C^{n-1}$ because both maps $J$ and $f \circ h$ are. Consequently the
composition map $\Phi \circ S_{J}^{f \circ h}: \Omega' \ra F$, mapping each $y \in
\Omega'$ into $(Jy) (f(h(y))) \in F$ is of class $C^{n-1}$. 

Now we check the form of  its first partial derivatives. We just see the partial derivative with respect to the first coordinate $x_1$. It is easy to check that
\begin{eqnarray*}
\frac{\partial}{\partial x_1} ( \Phi \circ S_{J}^{f \circ h} ) (y) &=&  \frac{\partial}{\partial X} \Phi ( S_{J}^{f \circ h} (y)) \circ \frac{\partial}{\partial x_1} Jy + 
 \frac{\partial}{\partial Y} \Phi (S_{J}^{f \circ h} (y)) \circ \frac{\partial}{\partial x_1} (f \circ h)( y)  \\ &=&
 \Phi \left( \frac{\partial}{\partial x_1}  Jy , f(h(y)) \right) + 
 \Phi \left( Jy,   \frac{\partial}{\partial x_1} ( f \circ h ) (  y ) \right) \\ &=& \left( \Phi \circ S_{\frac{\partial}{\partial x_1}  J}^{ f \circ h} \right) (y) + \left( \Phi \circ S_{J}^{\frac{\partial}{\partial x_1}  f \circ h} \right) (y).
\end{eqnarray*}

By an inductive reasoning, we see that the partial derivatives of order $n-1$ of $\Phi \circ S_{J}^{f \circ h}$ are just a sum of terms of the form $\Phi \circ S_{\partial^{\lambda} J}^{\partial^{\mu} (f \circ h)}$, where $\va \lambda \vb , \va \mu \vb \le n-1$. Consequently, to prove that $T f \in C^n (\Omega' ,F)$, we just have to show that each one of the terms $\Phi \circ S_{\partial^{\lambda} J}^{\partial^{\mu} (f \circ h)}$ ($\va \lambda \vb , \va \mu \vb \le n-1$) is of class $C^1$. We suppose that $S_K^{g} : \Omega' \ra L(E, F) \times E$ is one of the above $S_{\partial^{\lambda} J}^{\partial^{\mu} (f \circ h)}$ (that is, let us denote $g = \partial^{\mu} (f \circ h) \in C^1 (\Omega', E)$ and $K = \partial^{\lambda} J$), 
which is continuous when $L(E,F)$ is endowed with the topology of the norm, as we stated above.

Now we check the form of  its first partial derivatives. Take $i \in \{1, \ldots, p\}$, and  assume without loss of generality  that $i=1$. It is easy to check that
\[\lim_{k \ra 0} \frac{(\Phi \circ S_K^g) (y + (k, 0, \ldots, 0)) - (\Phi \circ S_K^g) (y)}{k} , \] that is, \[\lim_{k \ra 0} \frac{(K(y+ (k,0, \ldots, 0))) (g(y+ (k, 0 , \ldots, 0))) - (Ky) (g(y))}{k} \] is equal to 
\[ \lim_{k \ra 0} \frac{(K(y +(k, 0 , \ldots, 0))) (g(y+ (k, 0 , \ldots, 0))) - (K (y + (k, 0 , \ldots, 0)) )(g(y ))}{k} +\]\[+ \lim_{k \ra 0} \frac{(K ( y +(k, 0 , \ldots, 0)))(g(y )) - (Ky) (g(y))}{k} ,\]that is, it is equal to
\[(Ky) \left( \frac{\partial}{\partial x_1} g \left(y \right) \right) + \left( \frac{\partial_{\rm s}}{\partial x_1} K y \right) (g(y)) , \]by the definition of $\frac{\partial_{\rm s}}{\partial x_1} K y $ and the fact that $K$ is continuous for $L(E,F)$ endowed with the topology of the norm.

Applied to our context, we have that
\[\frac{\partial}{\partial x_1} ( \Phi \circ S_{\partial^{\lambda} J}^{\partial^{\mu} (f \circ h)} ) (y) =
 \left( \Phi \circ S_{\partial^{\lambda} J}^{\frac{\partial}{\partial x_1} \partial^{\mu} (f \circ h)} \right) (y) + \left( \Phi \circ S_{\frac{\partial_{\rm s}}{\partial_{ x_1} } \partial^{\lambda} J}^{ \partial^{\mu} (f \circ h)} \right) (y) .\]

Now, as we noted above, $S_{\partial^{\lambda} J}^{\frac{\partial}{\partial x_1} \partial^{\mu} (f \circ h)}$ is continuous when considering in $L(E,F)$ the topology of the norm. As a consequence, to obtain the continuity of all partial derivatives of order $n$ of $Tf$, it is enough to see that
$\Phi \circ S_{\frac{\partial_{\rm s}}{\partial_{ x_1} } \partial^{\lambda} J}^{ \partial^{\mu} (f \circ h)}$ is continuous.

In order to prove this, notice first that, since $J$ is of class ${\rm s}-C^{n-1}$, then for the above $\lambda$  the map 
$\frac{\partial_{\rm s}}{\partial_{ x_1} } \partial^{\lambda} J$ is continuous when $L(E, F)$ is endowed with the strong operator topology. This means that, given ${\bf e} \in E$, the map 
$\frac{\partial_{\rm s}}{\partial_{ x_1} } (\partial^{\lambda} J)_{\bf e} :\Omega' \ra F$ sending each $y \in \Omega'$ into $(\frac{\partial_{\rm s}}{\partial_{ x_1} } \partial^{\lambda} J y) ({\bf e})$ is continuous. Now take $\epsilon >0$ and $y_0 \in \Omega'$. We are going to show that there exists $\delta >0$ such that if $\va y-y_0 \vb < \delta$, $y \in \Omega'$, then
\[\vc \Phi \circ S_{\frac{\partial_{\rm s}}{\partial_{ x_1} } \partial^{\lambda} J}^{ \partial^{\mu} (f \circ h)} (y) - \Phi \circ S_{\frac{\partial_{\rm s}}{\partial_{ x_1} } \partial^{\lambda} J}^{ \partial^{\mu} (f \circ h)} (y_0) \vd < \epsilon.\]

Let ${\bf e}_0 := \partial^{\mu} (f \circ h) (y_0) $, and take $\delta_1 >0$
such that the closed ball $\bar{B} (y_0 , \delta_1)$ is contained in $\Omega'$. 
Since $\frac{\partial_{\rm s}}{\partial_{ x_1} } (\partial^{\lambda} J)_{{\bf e}_0}$ is continuous, there exists an upper bound $M$ for this function on the compact set $\bar{B} (y_0 , \delta_1)$. Also, there exists $\delta_2 >0$, $\delta_2 < \delta_1$, such that if $\va y - y_0 \vb < \delta_2$, then \[\vc \left( \frac{\partial_{\rm s}}{\partial_{ x_1} } \partial^{\lambda} J y 
-\frac{\partial_{\rm s}}{\partial_{ x_1} } \partial^{\lambda} J y_0
\right) \left( {\bf e}_0 \right) \vd < \frac{\epsilon}{2}.\]

On the other hand, since $\partial^{\mu} (f \circ h)$ is continuous, there exists $\delta >0$, $\delta < \delta_2 $, such that if $y \in B(y_0, \delta)$, then $\vc \partial^{\mu} (f \circ h) (y) - {\bf e}_0 \vd < \epsilon M/2$.  

Consequently, if $\va y - y_0 \vb < \delta$, we have
\[\vc \Phi \circ S_{\frac{\partial_{\rm s}}{\partial_{ x_1} } \partial^{\lambda} J}^{ \partial^{\mu} (f \circ h)} (y) - \Phi \circ S_{\frac{\partial_{\rm s}}{\partial_{ x_1} } \partial^{\lambda} J}^{ \partial^{\mu} (f \circ h)} (y_0) \vd \] 
is less than or equal to
\[\vc \left( \frac{\partial_{\rm s}}{\partial_{ x_1} } \partial^{\lambda} J y \right) \left( \partial^{\mu} (f \circ h) (y) -{\bf e}_0 \right)  \vd + 
\vc \left( \frac{\partial_{\rm s}}{\partial_{ x_1} } \partial^{\lambda} J y 
- \frac{\partial_{\rm s}}{\partial_{ x_1} } \partial^{\lambda} J y_0
\right) ({\bf e}_0)  \vd, \]which is easily  strictly less than $\epsilon$. This proves the continuity of our functions, as it was to see.
\end{pf}

The following result is a direct consequence of Lemma~\ref{401} and Theorem~\ref{401b}. Roughly speaking, it says that Proposition~\ref{dadada} provides the only way to construct linear biseparating maps from $\noe$ onto $\nmoe$. We state the result in its complete form.

\begin{theorem}\label{401c}
Suppose that $T: \noe \ra \moe$
is a ${\Bbb K}$-linear  biseparating map. Then $p=q$, $n=m$, and there
exist a diffeomorphism $h$ of class $C^n$ from $\Omega'$ onto
$\Omega$ and a  map
$J: \Omega' \ra L(E,F)$ of class ${\rm s}-C^n$ such that for every
$y \in \Omega'$ and 
every $f \in \noe$, \[ (Tf) (y) = (Jy) (f(h(y))) .\]Moreover, $Jy \in B(E, F)$ for every $y \in \Omega'$.
\end{theorem}

\begin{pf}
We will follow the same notation as in Lemma~\ref{401}, which provided a first description of linear biseparating maps, in particular everything related to the definition of $J$ and $h$. Also, by Lemma~\ref{401}, $p=q$ and $n=m$. 

First, for each $\lambda \in \Lambda$, we define a map $J_{\lambda}: \Omega' \ra L' (E, F)$ as  \[(J_{\lambda} y) ({\bf e}) :=  \partial^{\lambda} (T \widehat{\bf e}) (y),\] for $y \in \Omega'$ and ${\bf e} \in E$. $J_{\lambda}$ is clearly well defined.

The rest of the proof will apply just for the case when we are in Situation 1, but it is easy to see that slight changes in it allow to prove the theorem when we are in Situation 2. We will prove it through several claims.

\medskip

        {\bf Claim 1.} {\em For every $y \in \Omega'$ and every $\lambda \in  \Lambda$, $J_{\lambda} y$ belongs to $L(E,F)$.}

Take $y \in \Omega'$ and a sequence $({\bf e}_n)$ in $E$ converging to zero.  We will see that $((J_{\lambda} y) ({\bf e}_n) )$ goes to zero. First we have that, since $T$  is continuous by Theorem~\ref{401b}, the sequence of functions $(T \widehat{{\bf e}_n})$ converges to zero, which implies in particular that $((\partial^{\lambda} T \widehat{{\bf e}_n})(y))$ goes to zero. But this last sequence is precisely $((J_{\lambda} y) ({\bf e}_n))$, so the claim is proved.

\medskip

        {\bf Claim 2.} {\em For each $\lambda \in \Lambda$, $\va \lambda \vb \le n-1$, the map $J_{\lambda} : \Omega' \ra L(E,F)$ is continuous when  $L(E,F)$ is endowed with the topology of the norm.}

We will show that if $y_0 \in \Omega'$, then $J_{\lambda}$ is continuous at $y_0$. 
Since $T$ is continuous, we have that, for $r>0$ such that the closed ball $\bar{B} (y_0, r)$ is contained in $\Omega'$, there exists $M>0$ such that $p_{\bar{B}(y_0, r)} ( T \widehat{\bf e} ) < M$ holds for  every ${\bf e} \in E$ with $\vc {\bf e} \vd \le 1$. This implies that, for these ${\bf e}$, if $ \va y - y_0 \vb <r$, then $\vc (D \partial^{\lambda} T \widehat{\bf e})(y) \vd \le p M$. Consequently, as it can be seen for instance in \cite[Theorem 3.3.2]{Ca}, we have that  
\[  \vc (\partial^{\lambda} T\widehat{\bf e}) (y) - \partial^{\lambda} (T \widehat{\bf e}) (y_0) \vd \le p M \va y -y_0 \vb, \]for every ${\bf e} \in E$ with $\vc {\bf e} \vd \le 1$. Now, taking into account that $(J_{\lambda} y) ({\bf e}) = (\partial^{\lambda} T \widehat{\bf e}) (y)$ for every $y \in \Omega'$, the result follows, and the claim is proved.

\medskip

        {\bf Claim 3.} {\em For each $\lambda \in \Lambda$, $\va \lambda \vb \le n-1 $, $\partial^{\lambda} J = J_{\lambda}$.} 

Of course, the result is clear if $n=1$, so we suppose that $n \ge 2$.
We will just prove the claim in the particular case when $\lambda= \lambda_1: = (1, 0 , \ldots, 0)$. The proof for all other $\lambda \in \Lambda$ is similar and can be achieved inductively.

Take $y_0 \in \Omega'$ and $r>0$ such that the closed ball $\bar{B} (y_0, r) \subset \Omega'$. It is clear that if $h \in {\Bbb R}-\{0\}$, $\va h \vb <r$, and if ${\bf e}$ is in the closed unit ball of $E$, then  \[\vc \frac{J (y_0 + (h, 0, \ldots, 0) ) ({\bf e}) - (J y_0) ({\bf e})}{h}  -  (J_{\lambda_1} y_0) ({\bf e}) \vd \]
is equal to \[\vc \frac{(T \widehat{\bf e}) (y_0 + (h, 0, \ldots, 0)) - (T \widehat{\bf e}) (y_0)}{h} -  (\partial^{\lambda_1} T \widehat{\bf e}) (y_0) \vd, \] 
which, by \cite[Corollary XIII.4.4]{Ln}, is less than \[ \sup_{\va y - y_0\vb < 2h} \vc (\partial^{\lambda_1} T \widehat{\bf e}) (y_0) - (\partial^{\lambda_1} T \widehat{\bf e}) (y)  \vd,  \]that is, less than
\[ \sup_{\va y - y_0\vb < 2h} \vc J_{\lambda_1} y_0 - J_{\lambda_1} y  \vd. \]

Clearly this implies that
\[\lim_{h \ra 0} \vc \frac{ J (y+ (h, 0, \ldots, 0) ) - J y}{h} -  J_{\lambda_1} y \vd =0. \] Consequently, the partial derivative of $J$ with respect to the first coordinate exists at each $y_0 \in \Omega'$ and is equal to $J_{\lambda_1} y_0$.

\medskip

	{\bf Claim 4.} {\em Take $\lambda = (n_1, n_2 , \ldots, n_i , \ldots, n_p) \in \Lambda$ with $\va \lambda \vb = n-1$. Then $J_{\lambda}$ is of class ${\rm s}-C^1$. Moreover, if for each $i \in \{1, \ldots, n\}$, $\mu_i = (n_1 , n_2 , \ldots , n_i +1 , \ldots, n_p)$, then 
\[\frac{\partial_{\rm s}}{\partial x_i} J_{\lambda} = J_{\mu_i} .\]}

The proof that $\frac{\partial_{\rm s}}{\partial x_i} J_{\lambda} = J_{\mu_i}$ is similar to the proof of Claim 3 we have done, taking into account that $J_{\mu_i}$ is perhaps no longer continuous when $L(E,F)$ is endowed with its norm but $\partial^{\mu_i} T \widehat{\bf e}$ is continuous for every ${\bf e} \in E$. Consequently, to finish we just have to show that $J_{\mu}$, $\va \mu \vb =n$, is continuous when considering in $L(E,F)$ the strong operator topology. We have to prove that if $\mu \in \Lambda$, and $(y_n) $ is a sequence in $\Omega'$ converging to $y \in \Omega'$, then $(J_{\mu} y_n) ({\bf e})$ converges to $(J_{\mu} y) ({\bf e})$ for each ${\bf e} \in E$. But this is immediate from the definition of $J_{\mu}$.
\end{pf}

\noindent
{\bf Remark.} Notice that, in the case when we are in Situation 2, in Theorem~\ref{401c} the map $J$ and its partial (${\rm s}$-)derivatives up to  order $n$ can also be extended continuously to the boundary of $\Omega'$ in a natural way, when considering in $L(E,F)$ the strong operator topology.

\medskip

In the special case when $E={\Bbb K} =F$, we immediately deduce the following result.

\begin{corollary}
Suppose that $T: A^n (\Omega, {\Bbb K}) \ra A^m (\Omega', {\Bbb K})$
is a ${\Bbb K}$-linear  biseparating map. Then $p=q$, $n=m$, and there
exist a diffeomorphism $h$ of class $C^n$ from $\Omega'$ onto
$\Omega$ and a  map $a: \Omega' \ra {\Bbb K}$ of class $C^n$ which does not vanish at any point of $\Omega'$, such that for every
$y \in \Omega'$ and 
every $f \in A^m (\Omega, {\Bbb K})$, \[ (Tf) (y) = a(y) f(h(y)) .\]
\end{corollary}

\medskip

We finish with a corollary whose proof is easy from Theorem~\ref{401c}.

\begin{corollary}\label{poin}
If $\noe$ and $\nmoe$ are endowed with the topology of the pointwise convergence, then every linear biseparating map $T: \noe \ra \nmoe$ is continuous.
\end{corollary}

\bigskip

\noindent
{\bf Final Remark.}
Even if, in Sections 3, 4 and 5, we are considering two possible situations for our spaces of vector-valued functions, Lemma~\ref{401} and Theorems~\ref{401b} and~\ref{401c} can be given in a broader context. 
Our examples $\noe$ and $\moe$ give the pattern to follow to generalize the results to some other spaces. In particular, in Section 3, we analyzed the setting which ensured us the existence of the support map $h$. We had there two possible situations: in the first one, $\Omega$ and $\Omega'$ where open sets without any additional restrictions; in the second one, they were assumed to be also bounded. In this second case, we were assuming that functions in $\noe$ and $\moe$ admitted continuous extensions to the closure of $\Omega$ and $\Omega'$, respectively. According to these two following assumptions on $\Omega$ and $\Omega'$, sets $\Omega_1$ and $\Omega'_1$ where defined. We keep that notation in which follows, that is, $\Omega_1 := \Omega$ and $\Omega'_1 := \Omega'$ when they are not bounded, and $\Omega_1 := \bar{\Omega}$ and $\Omega'_1 :=\bar{\Omega'}$ when both  are bounded. This would lead us to the construction of a homeomorphism $h$ from $\Omega'_1$ onto $\Omega_1$ in some other cases. Namely, to construct the support map $h$, we need that subspaces $A (\Omega , E) \subset C^n (\Omega, E) \cap C(\Omega_1, E)$ and $B( \Omega', F) \subset C^m (\Omega', F) \cap C(\Omega'_1, F)$ are an $A$-module and a $B$-module, respectively, where $A \subset C^n (\Omega, {\Bbb R}) \cap C(\Omega_1, {\Bbb R})$ and $B \subset C^m (\Omega', {\Bbb R}) \cap C(\Omega'_1 , {\Bbb R})$ are strongly regular rings (see \cite{A1} for an appropriate description).

On the other hand, apart from these necessary conditions for constructing $h$, a careful reading of the proofs shows that the essential requirements that linear subspaces $A(\Omega, E) $ and $B(\Omega', F)$ must meet so as to satisfy Lemma~\ref{401} are: 1) $A(\Omega, E)$ and $B(\Omega, E)$ contain the constant functions; and 2) $C_c^n (\Omega, {\Bbb R}) \subset A$ and $C_c^m (\Omega', {\Bbb R}) \subset B$. In the case when $\Omega$ and $\Omega'$ are bounded and the coordinate projections $x_i$ belong to $A$ and $B$, then the proof of Lemma~\ref{401} can be followed step by step with no changes. Otherwise (even if $\Omega$ and $\Omega'$ are not bounded), the changes are few and natural.

As for Theorem~\ref{401b}, once we have that the above conditions are satisfied and our spaces can be endowed with a topology compatible with pointwise convergence, then a similar statement and proof holds for $A(\Omega, E)$ and $B(\Omega' , F)$ whenever $A(\Omega, E)$ contains all functions in $C^{n} (\Omega, E)$ with compact support.

Finally, if we want to obtain a  result similar to Theorem~\ref{401c} for our spaces $A(\Omega, E)$ and $B(\Omega', F)$, besides all the above conditions, they must  be endowed with a {\em suitable} norm or family of seminorms providing a topology compatible with the pointwise convergence. This will be the case, for instance, of the spaces of functions with bounded derivative.

\medskip

So we study the case of the spaces
$C^n_* (\Omega, E) \subset C^n (\Omega , E)$ and $C^m_* (\Omega', F) \subset C^m (\Omega' , F)$ consisting of all functions such that all partial derivatives up to orders $n$ and $m$, respectively, are {\em bounded}. The space $\bnoe$ (and similarly $\bmoe$)
becomes a Banach space with  the norm defined for each  $f \in \bnoe$ as \[\vc f\vd := \max_{\lambda
\in \Lambda} \sup_{x \in \Omega} \vc \partial^{\lambda}
f (x) \vd .\] This is a suitable norm in the above sense because the proof of Theorem~\ref{401c} can be followed easily for these spaces equipped with such norm.

First, notice that, since we are assuming that $n , m  \ge 0$, then in particular when $\Omega$ is convex all functions in $C^n_* (\Omega, E)$ admit a continuous extension to the closure of $\Omega $ in ${\Bbb R}^p$: suppose that $(x_n)$ is a sequence in $\Omega$ converging to $x_0$ in the boundary of $\Omega$, and that $f \in C^n_* (\Omega, E)$; then since the differential $D f$ is bounded on the whole $\Omega$ by an $M>0$, we have that, by
 \cite[Theorem 3.3.2]{Ca},
\[  \vc f (x_n) - f(x_m) \vd \le  M \va x_n - x_m \vb, \]which implies that  $(f(x_n))$ is a Cauchy sequence. In this way we would define the extension $f(x)$ as the limit of this sequence. It is straightforward to see that the new extended function is continuous in the closure of $\Omega$. 

On the other hand, when $\Omega$ and $\Omega'$ are {\em bounded} and {\em convex}, then it is easy to see that $C^n_* (\Omega, E)$ is a $C^n (\bar{\Omega}, {\Bbb R})$-module, and a similar statement is also valid for $C^m_* (\Omega', F)$. As a consequence, by the comments given above, Lemma~\ref{401} and Theorems~\ref{401b} and ~\ref{401c} can also be stated in this new situation (for $\Omega$ and $\Omega'$ bounded and convex). Furthermore, in this case, as in the remark after Theorem~\ref{401c}, it is also possible to say that partial derivatives up to order $n-1$ of $J$ admit a continuous extension (when $L(E,F)$ is equipped with the strong operator topology) to the boundary of $\Omega'$. As for the partial ${\rm s}$-derivatives of all partial derivatives of order $n-1$ of $J$, an elementary application of the Uniform Boundedness Theorem shows that they are bounded on $\Omega'$.

\medskip

What happens if $\Omega$ is for instance not bounded?
One might be tempted to follow a similar pattern as indicated above when trying to describe linear biseparating maps defined between $C^n_* (\Omega, E)$ and $C^m_* (\Omega', F)$.
But, in that case, we have that $\bnoe$ is no longer a $C^n (\Omega, {\Bbb R})$-module (as $\noe$ was). Anyway, it is a $C^n_* (\Omega, {\Bbb R})$-module, and we could try to follow the proof of Lemma~\ref{401} to get a similar description of linear biseparating maps, but even if we could manage to adapt the proof step by step (with some changes), there is a major problem from the beginning: in general the space $C^n_* (\Omega, {\Bbb R})$ is not a strongly regular ring, so our results cannot be applied, and in particular the existence of the support map $h$ is not clear. Let us see an example where $C^n_* (\Omega, {\Bbb R})$ is not a strongly regular ring.

\medskip

\noindent
{\em Example.} Suppose that $\Omega \subset {\Bbb R}$ is the open set defined as $\Omega := \bigcup_{k=1}^{\infty} (k - 1/k , k + 1/k)$. Let us define $K$ as the closure in $\beta \Omega$ of ${\Bbb N}$. Now take $x_0 \in \beta \Omega - (K \cup {\Bbb R})$. It is clear that if a derivable function $f: \Omega \ra {\Bbb R}$  satisfies $f^{\beta \Omega} \equiv 0$ on $K$ and $f^{\beta \Omega} \equiv 1$ on a neighborhood of $x_0$, then for some sequence $(x_n)$ in $\Omega$ going to infinity, the sequence $(f(x_n))$ converges to $1$. As a consequence from the Mean Value Theorem, we conclude that the derivative of $f$ cannot be bounded on $\Omega$, that is, $f \notin C^1_* (\Omega , {\Bbb R})$, as we wanted to show.

\medskip

We end the paper with some related questions, concerning special cases where our techniques cannot be applied.

\medskip

\noindent
{\bf Problem 1.} Assume that there exists a biseparating map $T: \noe \ra \moe$ which is not linear. Can we deduce that the support map $h: \Omega' \ra \Omega$ is a diffeomorphism of class $n$? Remark that the assumption of linearity in Lemma~\ref{401} is necessary for its proof.

\noindent
{\bf Problem 2.} Suppose that $C^{\infty} (\Omega, E)$ is the space of $E$-valued functions which are of class $C^{\infty}$ in $\Omega$, and that $C^{\infty} (\Omega', F)$ is defined in a similar way. Describe the  linear biseparating maps from $C^{\infty} (\Omega, E)$ onto $C^{\infty} (\Omega', F)$. Must such a map be continuous? Notice that by the comments given in the Final Remark above, the construction of the support map $h$ is possible, but the proof of Lemma~\ref{401} is no longer valid.

\noindent
{\bf Problem 3.} Let $\Omega$ and $\Omega'$ be unbounded open subsets of ${\Bbb R}^p$ and ${\Bbb R}^q$, respectively. Describe the  linear biseparating maps from $\bnoe$ onto $\bmoe$. 

\noindent
{\bf Problem 4.} Determine all subspaces $A(\Omega, E) \subset \noe$ and $B(\Omega' , F) \subset \moe$ such that the existence of a (linear) biseparating map from $A(\Omega , E)$ onto $B(\Omega', F)$ implies that $E$ and $F$ are isomorphic as Banach spaces.

\medskip

The author is thankful to Professor Richard Aron for drawing his attention to some papers which appear in the references.


\begin{thebibliography}{GHJRl}
\bibitem[AK]{AK} Y.\ A.\ Abramovich and A.\ K.\ Kitover, {\em Inverses of disjointness preserving operators}.\ Mem.\ Amer.\ Math.\ Soc.\ {\bf 143} (2000).
\bibitem[A1]{A1} J.\ Araujo, {\em Realcompactness and spaces of vector valued functions}. Submitted. Available at http://www.matesco.unican.es/\~{}araujo/homeo.dvi and at http://xxx.lanl.gov/pdf/math/0010261.
\bibitem[A2]{A2} J.\ Araujo, {\em Realcompactness and Banach-Stone theorems}. Submitted. Available at http://www.matesco.unican.es/\~{}araujo/isom3.dvi and at http://xxx.lanl.gov/pdf/math/0010292.
\bibitem[A3]{A3} J.\ Araujo, {\em Separating maps and linear
isometries between some spaces of continuous functions}.\ J.\
Math.\ Anal.\ Appl.\ {\bf 226} (1998), 23-39.
\bibitem[ABN]{ABN} J.\ Araujo, E.\ Beckenstein, and L.\ Narici,
{\em Biseparating maps and homeomorphic realcompactifications}.\
J.\ Math.\ Anal.\ Appl.\ {\bf 192} (1995), 258-265.
\bibitem[AGL]{AGL} R.\ Aron, J.\ G\'omez, and J.\ G.\ Llavona, {\em Homomorphisms between algebras of differentiable functions in infinite dimensions.}\ Michigan Math.\ J.\ {\bf 35} (1988), 163-178.
\bibitem[BCL]{BCL} W.\ G.\ Bade, P.\ C.\ Curtis, Jr., and K.\ B.\ Laursen, {\em  Automatic continuity in algebras of differentiable functions}.\
Math.\ Scand.\ {\bf 40} (1977), 249--270.
\bibitem[Ca]{Ca} H.\ Cartan, {\em Calcul diff\'erentiel}.\ Hermann,
Paris, 1967.
\bibitem[Ch]{Ch} G.\ Chilov, {\em Analyse Math\'ematique. Fonctions
 de plusieurs variables r\'eelles}, 2nd edition.\ Mir, Moscow, 1975.
\bibitem[D]{D} H.\ G.\ Dales, {\em Automatic continuity: a
survey}.\ Bull.\ London Math.\ Soc.\ {\bf 10} (1978), 129-183.
\bibitem[Fr]{Fr} H.\ Federer, {\em Geometric measure theory}.\ Springer Verlag, 1969.
\bibitem[Fk]{Fk} V.\ V.\ Fedorchuk, {\em The fundamentals of
dimension theory}.\ In {\em General Topology I}, edited by A.\ V.\ Arkhangelski\u{\i} and L.\ S.\ Pontryagin, Springer Verlag, 1990,
91-195.
\bibitem[Fl]{Fl} W.\ Fleming, {\em Functions of several
variables}, 2nd edition.\ Springer Verlag, 1977.
\bibitem[FH]{FH} J.\ J.\ Font and S.\ Hern\'andez, {\em On separating maps between locally compact spaces}.\ Arch.\ Math.\ (Basel) {\bf 63} (1994), no. 2, 158--165. 
\bibitem[GJ]{GJ} L.\ Gillman and M.\ Jerison, {\em Rings of
continuous functions}.\ Van Nostrand, Princeton, 1960.
\bibitem[GL]{GL} J.\ M.\ Guti\'errez and J.\ G.\  Llavona, {\em  Composition operators between algebras of differentiable functions}.\
Trans.\ Amer.\ Math.\ Soc.\ {\bf 338} (1993), 769--782.
\bibitem[J]{J} K.\ Jarosz, {\em Automatic continuity of
separating linear isomorphisms}.\ Canad.\ Math.\ Bull.\ {\bf 33}
(1990), 139-144.
\bibitem[JW]{JW}  J-S.\ Jeang and N-C.\ Wong, {\em Weighted composition operators of $C\sb 0(X)$'s}.\ J.\ Math.\ Anal.\ Appl.\ {\bf 201} (1996), 981--993. 
\bibitem[KN]{KN} R.\ Kantrowicz and  M.\ M.\ Neumann, {\em Automatic continuity of homomorphisms and derivations on algebras of vector-valued continuous functions}.\  Czecholovak  Math.\  J.\  {\bf 45} (1995), 747-756.
\bibitem[Ln]{Ln} S.\ Lang, {\em Real and Functional Analysis}, 3rd. 
edition.\ Springer-Verlag, New York, 1993.
\bibitem[Lr]{Lr} K.\ B.\ Laursen, {\em Prime ideals and automatic
continuity in algebras of differentiable functions}.\ J.\ Funct.\ Anal.\ 
{\bf 38} (1980), 16-24.
\bibitem[My]{My} S.\ B.\ Myers, {\em Algebras of differentiable functions}.\ Proc.\ Amer.\ Math.\ Soc.\ {\bf 5} (1954), 917-922.
\bibitem[Ml]{Ml} L.\ Molnar, {\em On rings of differentiable functions}.\ In {\em Contributions to the theory of functional equations, II}, 13--16. Grazer Math. Ber., {\bf 327}, Karl-Franzens-Univ. Graz, Graz, 1996.
\bibitem[NRV]{NRV} M. \ M.\ Neumann, A.\ Rodr\'{\i}guez-Palacios, and M.\ V.\ Velasco, {\em Continuity of homomorphisms and derivations on algebras of vector-valued functions}.\ Quart.\  J.\  Math.\ Oxford (2) {\bf 50} (1999), 279-300.
\bibitem[P]{P} L.\ E.\ Pursell, {\em Rings of continuous functions on open convex subsets of $R\sp{n}$}.\ Proc.\ Amer.\ Math.\ Soc.\ {\bf 19} (1968), 581--585.
\bibitem[S]{S} A.\ M.\ Sinclair, {\em Automatic continuity of
linear operators}.\ London Math.\ Soc.\ Lecture Notes Series 21,
Cambridge University Press, 1976.
\bibitem[W]{W} J.\ Wloka, {\em Partial differential equations}.\
Cambridge University Press, Cambridge-New York, 1987.
\end{thebibliography}
\end{document}